\documentclass[12pt]{article}
\usepackage{amsmath}
\usepackage{tikz}
\usepackage{float}
\usepackage{amssymb}
\usepackage{comment}
\usetikzlibrary{decorations.markings}
\tikzstyle{vertex}=[circle, draw, inner sep=0pt, minimum size=6pt]

\newtheorem{theorem}{Theorem}[section]

\newtheorem{lemma}[theorem]{Lemma}

\newtheorem{construction}[theorem]{Construction}

\newenvironment{pf}           {\noindent{\bf Proof.} }%
                                {\null\hfill$\Box$\par\medskip\medskip\medskip\medskip}

\begin{document}

%%%%%%%%%%%%%%%%%%%%%%%%%%%
%%  Article Begins       %%
%%%%%%%%%%%%%%%%%%%%%%%%%%%
\title{Applying Skolem Sequences to Gracefully Label New Families of Triangular Windmills}
\author{Ahmad H. Alkasasbeh \footnote{ahmad84@mun.ca} \hspace{0.5in} Danny Dyer \footnote{dyer@mun.ca} \hspace{0.5in} Nabil Shalaby \footnote{nshalaby@mun.ca}\\
Department of Mathematics and Statistics \\
Memorial University of Newfoundland \\
St. John's, Newfoundland \\
A1C 5S7 Canada}
\date{}
\maketitle
\begin{abstract}

 A function $f$ is a \textit{graceful labelling} of a graph $G=(V,E)$ with $m$ edges if $f$ is an injection $f:V\mapsto \{0,1,2,\dots,m\}$ such that each edge $uv \in E$ is assigned the label $|f(u)-f(v)|$, and no two edge labels are the same. If a graph G has a graceful labelling, we say that $G$ itself is graceful.

In this paper, we prove any Dutch windmill with three pendant triangles is (near) graceful, which settles Rosa’s conjecture for a new family of triangular cacti.
% we proved Rosa's conjecture for a new family of triangular cacti: Dutch windmills of any order with three pendant triangles.
\end{abstract}

{\bf Keywords:} Graceful Labellings; Skolem-type sequences.
%\clearpage
%%%%%%%%%%%%%%%%%%%%%%%%%
%   Intro               %
%%%%%%%%%%%%%%%%%%%%%%%%%
%\tableofcontents
%\newpage
%\chapter{Applying Skolem Sequences to Gracefully Label New Families of Triangular Windmills}\label{ch:2}
%{\let\thefootnote\relax\footnote{{Submitted for Discrete Applied Mathematics Journal.}}}
%%%%%%%%%%%%%%%%%%%%%%%%%%%%%%%%%%%%%%%%%%%%%%%%%%%Introduction%%%%%%%%%%%%%%%%%%%%%%%%%%%%%%%%%%%%%%%%%%%%%%%%%%%
\section{Introduction}

In 1963, Ringel \cite{ringel} posed the following problem, given an arbitrary tree $T$ with $m$ edges: Can $K_{2m+1}$ be decomposed into $2m+ 1$ copies of $T$?

In 1967, Rosa introduced a graph labelling called a $\beta$-{\it valuation}, and used the concept of $\beta$-valuations to aid attempts to prove Ringel's conjecture. Later, in 1972, Golomb renamed a $\beta$-valuation as a {\it graceful labelling}, as it is still known today.%In \cite{rosa1} Rosa introduced the definition of graceful labelling and proved that Ringel's conjecture is valid if every tree has a graceful labelling.

%\textbf{I need to talk here a little bit about why graceful labellings were interesting to Rosa, and any general conjectures about the idea of graceful labelling.}\parGraceful labellings of graphs are essential for solving real-life problems, such as finding configurations of antennas in radio-astronomy \cite{Bermond}.
In \cite{Bermond}, Bermond studied the graceful labelling of Dutch windmills. Dutch windmills are a subtype of triangular cacti, connected graphs whose blocks are all triangles.

In \cite{rosa}, Rosa stated the following conjecture:

%\hspace*{1em}
\begin{enumerate}
\item all triangular cacti with $n\equiv 0$,$1\ ($mod$\ 4)$ are graceful,
\item all triangular cacti with $n\equiv 2$,$3\ ($mod$\ 4)$ are nearly graceful.
\end{enumerate}

%\textbf{Here I need to dipping into triangular cacti, give some references to families of graphs that can be gracefully labelled, to show how broad it is, and mention Gallian's survey.}\par
%where a triangular snake,
In 1989, Moulton \cite{mou} proved Rosa's conjecture for a triangular snake, a type of triangular cactus whose block cutpoint graph is a path. Rosa, in \cite{rosa}, recommended using Skolem-type sequences to label various families of triangular cacti. In 2012, Dyer, Payne, Shalaby, and Wicks \cite{dyer} verified Rosa's conjecture for a new class of triangular cacti: Dutch windmills with at most two pendant triangles by using Skolem-type sequences as Rosa suggested. Gallian, in his survey (A Dynamic Survey of Graph Labeling) \cite{Gallian}, mentioned that a proof for all triangular cacti seems hopelessly difficult.

In this paper, we develop a system in which we add pendant triangles (in all possible configurations) to a family of gracefully labelled graphs to get new gracefully labelled families. We verify Rosa's conjecture for Dutch windmills of any order with three pendant triangles, which will offer an extension of the results in \cite{dyer}. We can conclude this work with the following theorem.
% Here we used Langford sequences to obtain Skolem and hooked Skolem sequences of a considerable size; this technique was introduced in \cite{dyer}.
%\pagebreak
%6666666666666666666666666666666666666666666666666666666666666666666666666666666666666666666666666666666666666666666666666666666
\begin{theorem}
Every Dutch windmill with at most three pendant triangles is graceful or near graceful.
\end{theorem}
%666666666666666666666666666666666666666666666666666666666666666666666666666666666666666666666666666666666666666666666nnnnnnnnnnnnnnnnnnnnnnnn6666666666
%66666666666666666666666666666666666666666666666666666666666666666666666666666666666666666666666666666666666666666nnnnnnnnnnnnnnnnnnnnnnnn66666666666666
%66666666666666666666666666666666666666666666666666666666666666666666666666666666666666666666666666666666666nnnnnnnnnnnnnnnnnnnnnnnn66666666666666666666
\section{Definitions and Preliminaries}
Essential definitions are introduced in this section, as well as previously found results required to prove further results.
\subsection{Skolem and Langford Sequences}\label{skolemintro}
%\noindent In this paper, we use the definitions from the Handbook of Combinatorial Designs \cite{crc}.
The following definitions come from the Handbook of Combinatorial Designs \cite{crc}.

\rm A {\it Skolem sequence} of order $n$ is a sequence $S=\left(s_{1},s_{2},\ldots,s_{2n}\right)$ of $2n$ integers satisfying these conditions:

\begin{enumerate}

\item for every $k\in\left\{1,2,\ldots,n\right\},$ there exist exactly two elements $s_{i},s_{j}\in S$ such that $s_{i}=s_{j}=k;$

\item if $s_{i}=s_{j}=k,$ with $i < j,$ then $j-i=k$.

\end{enumerate}

For example, $S_4=(1,1,4,2,3,2,4,3)$ or equivalently $\left\{\left(1,2\right)\right.$, $\left(4,6\right)$, $\left(5,8\right)$, $\left.\left(3,7\right)\right\}$ is a Skolem sequence of order $4$. If we have a Skolem sequence $\{(a_i,b_i)\}_{i=1}^n,$ then $i$ is called a \textit{pivot} of a Skolem sequence if $b_i+i \leq 2n.$
\begin{theorem} \cite{skolem}
A Skolem sequence of order $n$ exists if and only if $n\equiv0,1\ ({\rm mod}\ 4)$.
\end{theorem}
{\rm A {\it hooked Skolem sequence} of order $n$ is a sequence $hS=\left(s_{1},s_{2},\ldots,s_{2n+1}\right)$ of $2n+1$ integers satisfying these conditions:

\begin{enumerate}

\item for every $k\in\left\{1,2,\ldots,n\right\},$ there exist exactly two elements $s_{i},s_{j}\in hS$ such that $s_{i}=s_{j}=k;$

\item if $s_{i}=s_{j}=k,$ with $i < j,$ then $j-i=k;$

\item $s_{2n}=0.$

\end{enumerate}}

For example, $hS_2=(1,1,2,0,2)$ is a hooked Skolem sequence of order $2$. If we have a hooked Skolem sequence $\{(a_i,b_i)\}_{i=1}^n,$ then $i$ is called a \textit{pivot} of a hooked Skolem sequence if $2n \neq b_i+i \leq 2n+1.$
\begin{theorem}\cite{okeefe}
A hooked Skolem sequence of order $n$ exists if and only if $n\equiv2,3\ ({\rm mod}\ 4)$.
\end{theorem}
%%%%%%%%%%%%%%%%%%%%%%%%%%%%%%%%%%%%%%   pivot         %%%%%%%%%%%%%%%%%%%%%%%%%%%%%%%%%%%%%%%%%%%%%%%%%%%%%%%%%%%%%
%If we have a (hooked) Skolem sequence $\{(a_i,b_i)\}_{i=1}^n,$ then $i$ is called a \textit{pivot} of a Skolem sequence if $b_i+i \leq 2n,$ and $i$ is called a \textit{pivot} of a hooked Skolem sequence if $2n \neq b_i+i \leq 2n+1.$\par

%\begin{enumerate}
%Let $i$ be a \textit{pivot} of a Skolem sequence $\{(a_i,b_i)\}_{i=1}^n$ if $b_i+i \leq 2n;$
%Let $i$ be a \textit{pivot} of a hooked Skolem sequence $\{(a_i,b_i)\}_{i=1}^n$ if $2n \neq b_i+i \leq 2n+1.$
%\end{enumerate}
%Assume we have a (hooked) Skolem sequence Let $\{(a_i,b_i)\}_{i=1}^n$  $i$ is said to be a \textit{pivot} of a Skolem sequence if $b_i+i \leq 2n$. For $i$ to be a pivot of a hooked Skolem sequence, $2n \neq b_i+i \leq 2n+1$.
%%%%%%%%%%%%%%%%%%%%%%%%%%%%%%%%%%%%%%%%%%%% Heffters %%%%%%%%%%%%%%%%%%%%%%%%%%%%%%%%%%%%%%%%%%%%%%
%\textbf{I need to discuss it with my supervisor ??? see CRC book page 614 and Ref \cite{dyer}}
In $1897$, Heffter stated two difference problems \cite{heffter1}. Heffter's first difference problem is: can a set $\{1,\ldots, 3n\}$ be partitioned into $n$ ordered triples $(a_i,b_i,c_i)$, with $1\leq i\leq n$, such that $a_i+b_i=c_i$ or $a_i+b_i+c_i\equiv 0\;(\textup{mod}\;6n+1)$? If such a partition is possible, then $\{\{0,a_i+n,b_i+n\}|1\leq i\leq n\}$ will be the base blocks of a cyclic Steiner triple system of order $6n+1$, $CSTS(6n+1)$. Construction \ref{ds} gives a solution to Heffter's first difference problem.
%\rm {Heffter's second difference problem} is: can a set $\{1,\ldots, 3n+1\}\backslash \{2n+1\}$  be partitioned into $n$ ordered triples $(a_i,b_i,c_i)$ with $1\leq i\leq n$ such that $a_i+b_i=c_i$ or $a_i+b_i+c_i\equiv 0\ (\textup{mod}\;6n+3)$? If such a partition is possible then $\{\{0,a_i+n,b_i+n\}|1\leq i\leq n\}$ with the addition of the base block $\{0,2n+1,4n+2\}$ having a short orbit of length $3n+1$ will be the base blocks of a $CSTS(6n+3)$.
\begin{construction}\cite{skolem}{\label{ds}}
Consider the (hooked) Skolem sequence with pairs $(a_i,b_i)$. The set of all triples $(i,a_i+n,b_i+n)$, for $1\leq i\leq n$, is a solution to the Heffter first difference problem. These triples yield the base blocks for a $CSTS(6n+1)$: $\{0,a_i+n,b_i+n\}$, $1\leq i\leq n$. Also, $\{0,i,b_i+n\}$, $1\leq i\leq n$ is another set of base blocks of a $CSTS(6n+1)$.
\end{construction}

{\label{baseblocks}}Let $S_4=(1,1,4,2,3,2,4,3)$ be a Skolem sequence of order $4$, yielding the pairs $\left(1, 2\right)$, $\left(4, 6\right)$, $\left(5, 8\right)$, $\left(3, 7\right)$. These pairs yield in turn the triples $\left(1, 5, 6\right)$, $\left(2, 8, 10\right)$, $\left(3, 9, 12\right)$, $\left(4, 7, 11\right)$, forming a solution to the first Heffter problem.
These triples yield the base blocks for two $CSTS(25)s$:
\begin{enumerate}
\item $\{0,5,6\},\{0,8,10\},\{0,9,12\},$ and $\{0,7,11\}\ (\textup{mod}\;25)$;
\item $\{0,1,6\},\{0,2,10\},\{0,3,12\},$ and $\{0,4,11\}\ (\textup{mod}\;25)$.
\end{enumerate}

%\textbf{I think we need to add the definition of Rosa and hooked Rosa ???}
%\\
%\begin{construction}\cite{rosa2}{\label {dhs}}
%A cyclic Steiner triple system $CSTS(6n+3),n> 1$. From a Rosa sequence or a hooked Rosa sequence of order $n$, construct the pairs $(a_i,b_i)$ such that $b_i-a_i=i$ for $1\leq i\leq n$. The set of all triples $(i,a_i+n,b_i+n)$ for $1\leq i\leq n$ is a solution to the Heffter second difference problem. These triples yield the base blocks for a $CSTS(6n+3)$: $\{0,a_i+n,b_i+n\}$, $1\leq i\leq n$ together with the base block $\{0,2n+1,4n+2\}$ having short orbit of length $3n+1$.\end{construction}
%%%%%%%%%%%%%%%%%%%%%%%%%%%%%%%%%%%%%%%%%%%%%%%%%%%%%%%%%%%%%%%%%%%%%%%%%%%%%%%%%%%%%%%%%%%%%%%%%%%%
%\subsection{Langford Sequences}
A \textit{Langford sequence} of defect $d$ and order $l$ is a sequence $L=\left(l_{1}\right.$, $l_{2}$, $\ldots$, $\left.l_{2l}\right)$ which satisfies these conditions:
\begin{enumerate}
	\item for every $k\in\left\{d,d+1,\ldots,d+l-1\right\},$ there exist exactly two elements $l_{i},l_{j}\in S$ such that $l_{i}=l_{j}=k$;
\item if $l_{i}=l_{j}=k,$ with $i < j,$ then $j-i=k$.
\end{enumerate}
%\end{dfn}
%\begin{dfn}\cite{Bermond}
A \textit{hooked Langford sequence} of defect $d$ and order $l$ is a sequence $L=\left(l_{1}, l_{2}\right.$, $\cdots$, $\left.l_{2l+1}\right)$ which satisfies these conditions:
\begin{enumerate}
	\item for every $k\in\left\{d,d+1,\ldots,d+l-1\right\},$ there exist exactly two elements $l_{i},l_{j}\in S$ such that $l_{i}=l_{j}=k$;

\item if $l_{i}=l_{j}=k,$ with $i < j,$ then $j-i=k;$
\item $l_{2m}=0.$
% $hL_{3}^{6} = (8,4,7,3,6,4,3,5,8,7,6,0,5).$
\end{enumerate}
For example, $(4,2,3,2,4,3)$ is a Langford sequence with $d=2$ and $l=3$ and $(8,4,7,3,6,4,3,5,8,7,6,0,5)$
 is a hooked Langford sequence with $d=3$ and $l=6$.
%\end{dfn}

%\begin{exm} $L=\left(7,5,3,6,4,3,5,7,4,6\right)$ is a Langford sequence of order $5$ and defect $3$.
%\end{exm}
The necessary and sufficient conditions for the existence of (hooked) Langford sequences are given in Theorem~\ref{langfordexistence}.
\begin{theorem}\cite{simpson}\label{langfordexistence} \hspace{1cm}
\begin{enumerate}
\item A Langford sequence of order $m$ and defect $d$ exists if and only if
	\begin{enumerate}
		\item $m \geq 2d-1$
		\item $m \equiv 0,1 \pmod 4$ and $d$ is odd, or
		\item $m \equiv 0,3 \pmod 4$ and $d$ is even.
	\end{enumerate}
\item A hooked Langford sequence of order $m$ and defect $d$ exists if and only if
	\begin{enumerate}
		\item $m(m-2d+1)+2 \geq 0$
		\item $m \equiv 2,3 \pmod 4$ and $d$ is odd, or
		\item $m \equiv 1,2 \pmod 4$ and $d$ is even.
	\end{enumerate}
\end{enumerate}
\end{theorem}
Though they can be thought of as a natural generalization of Skolem sequences, Langford sequences have also been classically used to build new Skolem sequences by concatenating sequences of appropriate order, possibly interlacing hooks, if required. This classic technique is a much used method to construct Skolem sequences from Langford sequences.
\begin{lemma}\label{skla} If a (hooked) Skolem sequence of order $d-1$ exists, and a (hooked) Langford sequence of order $l$ and defect $d,$ then a (hooked) Skolem sequence of order $N=l+d-1$ exists. In particular, a new Skolem sequence of order $N$ is obtained by concatenating a Skolem sequence with a Langford sequence or by interlacing a hooked Skolem sequence and hooked Langford sequence. A new hooked Skolem sequence of order $N$ is obtained by concatenating a Skolem sequence with a hooked Langford sequence or a hooked Skolem sequence and a Langford sequence. \end{lemma}
%\begin{proof}
%We present a proof for the first case: that is, if we have a Skolem sequence and a Langford sequence then we can obtain a Skolem sequence. The proof of the other three cases will be similar.\par
%Let $S_{d-1}$ be a Skolem sequence of order $d-1$ containing the numbers $(1,2,\ldots,d-1)$ and denoted $S_{d-1}.$ Theorem~\ref{langfordexistence} tells us that we can find a Langford sequence $L$ of length $m$ and defect $d$ containing the numbers $(d,d+1,\ldots,m+d-1)$. When $d \equiv 1 \text{ or } 2 \pmod 4,$ we can construct a new Skolem sequence of order $d-1 \equiv 0 \text{ or } 1 \pmod 4.$ Simply by concatenating $S_{d-1}$ to $L$ we obtain a new Skolem sequence of order $N=m+d-1.$
%\end{proof}
Let $hS_{2}=(1, 1, 2, 0, 2)$ be a hooked Skolem sequence and $hL_{3}^{6} = (8$, 4, 7, 3, 6, 4, 3, 5, 8, 7, 6, 0, $5)$ be a hooked Langford sequence. Now, if we take the reverse of the hooked Skolem sequence and by interlacing $hL_{3}^{6}$ and the reverse of $hS_{2},$ then we obtain a new Skolem sequence $S_{8}=(8$, 4, 7, 3, 6, 4, 3, 5, 8, 7, 6, 2, 5, 2, 1, $1)$ of order $8.$
%\subsection{Constructing Skolem Sequences from Langford Sequences}
%7777777777777777777777777777777777777777777777777777777777777777777777777777777777777777777777777777777777777777777777777777777

\subsection{Graceful Labellings and Triangular Cacti}
Let $G=(V,E)$ be a graph with $m$ edges. Let $f$ be a labelling defined from $V(G)$ to $\{0,1,2,\dots,m\}$ and let $g$ be the induced edge labelling defined from $E(G)$ to $\{1,2,\dots,m\}$ by $g(uv)=|f(u)-f(v)|$, for all $uv \in E.$ The labelling $f$ is said to be \textit{graceful}, if $f$ is an injective mapping and $g$ is a bijection. If a graph $G$ has a graceful labelling, then it is graceful.

Let $G=(V,E)$ be a graph with $m$ edges. Let $f$ be a labelling defined from $V(G)$ to $\{0,1,2,\dots,m+1\}$ and let $g$ be the induced edge labelling defined from $E(G)$ to $\{1,2,\dots,m-1,m\}$ or $\{1,2,\dots,m-1,m+1\}$ by $g(uv)=|f(u)-f(v)|$, for all $uv \in E.$ The labelling $f$ is said to be \textit{near graceful}, if $f$ is an injective mapping and $g$ is a bijection. If a graph $G$ has a near graceful labelling, then it is near graceful.
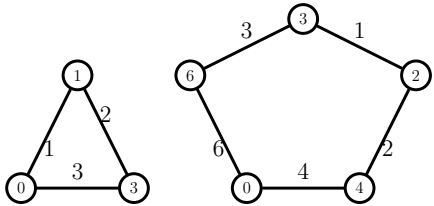
\begin{figure}[h!]
\begin{center}
\scalebox{0.75}{
\begin{tikzpicture}
     \node[shape=circle,draw=black, scale=0.7,line width=1.5pt] (0) at (0,0) { 0};
    \node[shape=circle,draw=black, scale=0.7,line width=1.5pt] (3) at (2,0) { 3};
   		   \node[shape=circle,draw=black, scale=0.7,line width=1.5pt] (1) at (1,2) { 1};	
		      \path [line width=1.5pt](0) edge[black] node[right,above] {\color{black} $3$} (3);
		      \path [line width=1.5pt](0) edge[black] node[left,below] {\color{black} $1$} (1);
		      \path [line width=1.5pt](1) edge[black] node[left,above] {\color{black} $2$} (3);
		 \node[shape=circle,draw=black, scale=0.7,line width=1.5pt] (0) at (4,0)  {\color{black} $0$};
    \node[shape=circle,draw=black, scale=0.7,line width=1.5pt] (4) at (6,0)   {\color{black} $4$};
		 \node[shape=circle,draw=black, scale=0.7,line width=1.5pt] (6) at (3,2)  {\color{black} $6$};
   \node[shape=circle,draw=black, scale=0.7,line width=1.5pt] (2) at (7,2)   {\color{black} $2$};
		 \node[shape=circle,draw=black, scale=0.7,line width=1.5pt](3) at (5,3)   {\color{black} $3$};	
		      \path [line width=1.5pt](0) edge[black] node[right,below] {\color{black} $6$} (6);
		      \path [line width=1.5pt](0) edge[black] node[left,above] {\color{black} $4$} (4);
		      \path [line width=1.5pt](4) edge[black] node[left,below] {\color{black} $2$} (2);
		      \path [line width=1.5pt](2) edge[black] node[left,above] {\color{black} $1$} (3);
		      \path [line width=1.5pt](3) edge[black] node[left,above] {\color{black} $3$} (6);
\end{tikzpicture}
}
\caption{Graceful labelling of $K_{3}$ and near graceful labelling of $C_{5}$.}\label{c}
\end{center}
\end{figure}
%Let $G=(V,E)$ be a graph with $m$ edges. Let $f$ be a labelling defined from $V(G)$ to $\{0,1,2,\dots,m-1\} \cup \{\ m\  $or$\  m+1\}$ ($f:V\mapsto \{0,1,2,\dots,m-1\} \cup \{\ m\  $or$\  m+1\}$) and let $g$ be the induced edge labelling defined from $E(G)$ to $\{1,2,\dots,m-1\} \cup \{\ m\  $or$\  m+1\}$ by $g(uv)=|f(u)-f(v)|_{{uv} \in E}.$ The labelling $f$ is said to be \textit{almost graceful}, if $f$ is a $1-1$ mapping and $g$ is a bijection. If a graph $G$ has a graceful labelling, then it is almost graceful.\par
In this paper, we will gracefully label some new families of triangular cacti and introduce some definitions and results related to triangular cacti.

{\rm A {\it triangular cactus} is a connected graph whose blocks are all triangles ($K_3$)}. {\rm A triangular cactus that has the property of all its blocks having a common vertex is said to be a \textit{Dutch windmill}, and the blocks will be called \textit{vanes}.} A {\it pendant triangle} is defined as a block that is added to any triangular cacti. Sometimes in this paper we call a Dutch windmill has $k$ vanes by $k$-vane Dutch windmill.

The necessary conditions for gracefulness and near gracefulness of general triangular cacti are as follows:
\begin{theorem}\cite{dyer}\label{necessarycacti}
Let $G$ be a triangular cactus with $n$ blocks. Then,
\begin{enumerate}
\item if $G$ is graceful, then $n \equiv 0,1 \pmod 4$, and,
\item if $G$ is near graceful, then $n \equiv 2,3 \pmod 4$.
\end{enumerate}\end{theorem}
%\textbf{I need to add three simple graphs as what my supervisor told me before}\\
%^^^^^^^^^^^^^^^^^^^^^^^^^^^^^^^^^^^^^^^^^^^^^^^^^^^^^^^^^^^^^^^^^^^^^^^^^^^^^^^^^^^^^^^^^^^^^^^^^^^^^^^^^^^^^^^^^^^^^^^^^^^^^^^^^^^^^^^^^^^^^^^^^^^^^^^^^^^^^^^^^^^^^^^^^^^^^^^^^^^^^^^^^^^^^^^^^^^^^^^^^^^^^^^^^^^^^^^^^^^^^^^^^^^^^^^^^^^^^^^^^^^^^^^^^^^^^^^^^^^^^^^^^^^^^^^^^^^^^^^^^^^^^^^^^^^^^^^^^^^^^^^^^^^^^^^^^^^^^^^^^^^^^^^^^^^^^^^^^^^^^^^^^^^^^^^^^^^^^^^^^^^^^^^^^^^^^^^^^^^^^^^^^^^^^^^^^^^^^^^^^^^^^^^^^^^^^^^^^^^^^^^^^^^^^^^^^^^^^^^^^^^^^^^^^^^^^^^^^^^^^^^^^^^^^^^^^^^^^^^^^^^^^^^^^^^^^^^^^^^^^^^^^^^^^^^^^^^^^^^^^^^^^^^^^^^^^^^^^^^^^^^^^^^^^^^^^^^^^^^^^^^^^^^^^^^^^^^^^^^^^^^^^^^^^^^^^^^^^^^^^^^^^^^
%\newpage
\section{Dutch Windmills with Three Pendant Triangles}\label{sectiondutch}
In order to verify Rosa's conjecture for a new family of triangular cacti, namely Dutch windmills of any order with three pendant triangles, we will use Langford sequences to obtain Skolem and hooked Skolem sequences of considerable sizes. This technique was introduced in \cite{dyer}. We categorize all such cacti into one of eleven types, called Type $(a)$ through Type $(k)$, and then gracefully label each type.
\begin{figure}[H]
\begin{center}
\begin{picture}(450,450)
\linethickness{1pt}
%\put(0,0){\line(1,0){260}}
%\put(0,0){\line(0,1){260}}
%\put(0,260){\line(1,0){260}}
%\put(260,0){\line(0,1){260}}
%\put(130,0){\line(0,1){260}}
%\put(0,130){\line(1,0){260}}
\put(25,335){$(a)$}
\put(220,335){$(b)$}
\put(25,205){$(c)$}
\put(220,205){$(d)$}
\put(25,85){$(e)$}   %
\put(220,85){$(f)$}

%\tikzstyle{every node}=[draw,shape=circle];
%%%%%%%%%%%%%%%%%%%%%%%%%%%%%%%%%%%%% Graph a %%%%%%%%%%%%%%%%%%%%%%%%%%%%%%%%%%%%%%%%%%%%%%%%%%%%%%%%%%%%%%%%%%%%%%%%%%%%%%%%%%%%%
%\put(25,630){\begin{tikzpicture}[scale=0.9]
%\tikzstyle{every node}=[draw,shape=circle];
\put(25,345){\begin{tikzpicture}[scale=1.1]
\path (0:0cm) node[draw,shape=circle] (v0) {};
\path (0:1cm) node[draw,shape=circle] (v1) {};
\path (40:1cm) node[draw,shape=circle] (v2) {};
\path (92:1cm) node[draw,shape=circle] (v3) {};
\path (132:1cm) node[draw,shape=circle] (v4) {};
\path (174:1cm) node[draw,shape=circle] (v5) {};
\path (214:1cm) node[draw,shape=circle] (v6) {};
\path (266:1cm) node[draw,shape=circle] (v7) {};
\path (306:1cm) node[draw,shape=circle] (v8) {};
%%%%%%%%%%%%%%%%%%%%%%%%%%%%%%%%%%%%%%%%%%%%%%%%%%%%%%%%%%%%%%%%%%%%%%%%%%%%%%%%
\path (85:1.931cm) node[draw=blue,shape=circle] (v9) {};
\path (110:1.931cm) node[draw=blue,shape=circle] (v10){};

\path (160:1.931cm) node[draw=red,shape=circle] (v11) {};
\path (180:2cm) node[draw=red,shape=circle] (v12) {};

\path (15:1.931cm) node[draw=green,shape=circle] (v13) {};
\path (30:1.732cm) node[draw=green,shape=circle] (v14) {};
\draw
(v0) -- (v1)
(v1) -- (v2)
(v2) -- (v0)
(v0) -- (v3)
(v3) -- (v4)
(v4) -- (v0)
(v0) -- (v5)
(v5) -- (v6)
(v6) -- (v0)
(v0) -- (v7)
(v7) -- (v8)
(v8) -- (v0)

(v5) edge[red] (v11)
(v5) edge[red] (v12)
(v11) edge[red] (v12)

(v1) edge[green] (v13)
(v1) edge[green] (v14)
(v13) edge[green] (v14)
%%%%%%%%%%%%%%%%%%%%
(v9) edge[blue] (v3)
(v10) edge[blue] (v3)
(v9) edge[blue] (v10);
\end{tikzpicture}}

%%%%%%%%%%%%%%%%%%%%%%%%%%%%%%%%%%%%% Graph b %%%%%%%%%%%%%%%%%%%%%%%%%%%%%%%%%%%%%%%%%%%%%%%%%%%%%%%%%%%%%%%%%%%%%%%%%%%%%%%%%%%%%
\put(220,345){\begin{tikzpicture}[scale=1.1]
%\put(220,630){\begin{tikzpicture}[scale=0.9]
%\tikzstyle{every node}=[draw,shape=circle];
\path (0:0cm) node[draw,shape=circle] (v0) {};
\path (0:1cm) node[draw,shape=circle] (v1) {};
\path (40:1cm) node[draw,shape=circle] (v2) {};
\path (92:1cm) node[draw,shape=circle] (v3) {};
\path (132:1cm) node[draw,shape=circle] (v4) {};
\path (174:1cm) node[draw,shape=circle] (v5) {};
\path (214:1cm) node[draw,shape=circle] (v6) {};
\path (266:1cm) node[draw,shape=circle] (v7) {};
\path (306:1cm) node[draw,shape=circle] (v8) {};
%%%%%%%%%%%%%%%%%%%%%%%%%%%%%%%%%%%%%%%%%%%%%%%%%%%%%%%%%%%%%%%%%%%%%%%%%%%%%%%%
\path (45:1.931cm) node[draw=blue,shape=circle] (v9) {};
\path (60:1.931cm) node[draw=blue,shape=circle] (v10){};

\path (160:1.931cm) node[draw=red,shape=circle] (v11) {};
\path (180:2cm) node[draw=red,shape=circle] (v12) {};

\path (0:1.931cm) node[draw=green,shape=circle] (v13) {};
\path (-15:1.732cm) node[draw=green,shape=circle] (v14) {};
\draw
(v0) -- (v1)
(v1) -- (v2)
(v2) -- (v0)
(v0) -- (v3)
(v3) -- (v4)
(v4) -- (v0)
(v0) -- (v5)
(v5) -- (v6)
(v6) -- (v0)
(v0) -- (v7)
(v7) -- (v8)
(v8) -- (v0)

(v5) edge[red] (v11)
(v5) edge[red] (v12)
(v11) edge[red] (v12)

(v1) edge[green] (v13)
(v1) edge[green] (v14)
(v13) edge[green] (v14)
%%%%%%%%%%%%%%%%%%%%
(v9) edge[blue] (v2)
(v10) edge[blue] (v2)
(v9) edge[blue] (v10);
\end{tikzpicture}}

%%%%%%%%%%%%%%%%%%%%%%%%%%%%%%%%%%%%% Graph c %%%%%%%%%%%%%%%%%%%%%%%%%%%%%%%%%%%%%%%%%%%%%%%%%%%%%%%%%%%%%%%%%%%%%%%%%%%%%%%%%%%%%
\put(25,215){\begin{tikzpicture}[scale=1.1]
%\put(25,510){\begin{tikzpicture}[scale=0.9]
%\tikzstyle{every node}=[draw,shape=circle];
\path (0:0cm) node[draw,shape=circle] (v0) {};
\path (0:1cm) node[draw,shape=circle] (v1) {};
\path (40:1cm) node[draw,shape=circle] (v2) {};
\path (92:1cm) node[draw,shape=circle] (v3) {};
\path (132:1cm) node[draw,shape=circle] (v4) {};
\path (174:1cm) node[draw,shape=circle] (v5) {};
\path (214:1cm) node[draw,shape=circle] (v6) {};
\path (266:1cm) node[draw,shape=circle] (v7) {};
\path (306:1cm) node[draw,shape=circle] (v8) {};
%%%%%%%%%%%%%%%%%%%%%%%%%%%%%%%%%%%%%%%%%%%%%%%%%%%%%%%%%%%%%%%%%%%%%%%%%%%%%%%%
\path (85:1.931cm) node[draw=blue,shape=circle] (v9) {};
\path (110:1.931cm) node[draw=blue,shape=circle] (v10){};

\path (345:1.931cm) node[draw=red,shape=circle] (v11) {};
\path (0:2cm) node[draw=red,shape=circle] (v12) {};

\path (15:1.931cm) node[draw=green,shape=circle] (v13) {};
\path (30:1.732cm) node[draw=green,shape=circle] (v14) {};
\draw
(v0) -- (v1)
(v1) -- (v2)
(v2) -- (v0)
(v0) -- (v3)
(v3) -- (v4)
(v4) -- (v0)
(v0) -- (v5)
(v5) -- (v6)
(v6) -- (v0)
(v0) -- (v7)
(v7) -- (v8)
(v8) -- (v0)

(v1) edge[red] (v11)
(v1) edge[red] (v12)
(v11) edge[red] (v12)

(v1) edge[green] (v13)
(v1) edge[green] (v14)
(v13) edge[green] (v14)
%%%%%%%%%%%%%%%%%%%%
(v9) edge[blue] (v3)
(v10) edge[blue] (v3)
(v9) edge[blue] (v10);
\end{tikzpicture}}

%%%%%%%%%%%%%%%%%%%%%%%%%%%%%%%%%%%%% Graph d %%%%%%%%%%%%%%%%%%%%%%%%%%%%%%%%%%%%%%%%%%%%%%%%%%%%%%%%%%%%%%%%%%%%%%%%%%%%%%%%%%%%%

\put(220,215){\begin{tikzpicture}[scale=1.1]
%\put(220,510){\begin{tikzpicture}[scale=0.9]
%\tikzstyle{every node}=[draw,shape=circle];
\path (0:0cm) node[draw,shape=circle] (v0) {};
\path (0:1cm) node[draw,shape=circle] (v1) {};
\path (40:1cm) node[draw,shape=circle] (v2) {};
\path (92:1cm) node[draw,shape=circle] (v3) {};
\path (132:1cm) node[draw,shape=circle] (v4) {};
\path (174:1cm) node[draw,shape=circle] (v5) {};
\path (214:1cm) node[draw,shape=circle] (v6) {};
\path (266:1cm) node[draw,shape=circle] (v7) {};
\path (306:1cm) node[draw,shape=circle] (v8) {};
%%%%%%%%%%%%%%%%%%%%%%%%%%%%%%%%%%%%%%%%%%%%%%%%%%%%%%%%%%%%%%%%%%%%%%%%%%%%%%%%
\path (30:2.631cm) node[draw=blue,shape=circle] (v9) {};
\path (40:2.332cm) node[draw=blue,shape=circle] (v10) {};

\path (160:1.931cm) node[draw=red,shape=circle] (v11) {};
\path (180:2cm) node[draw=red,shape=circle] (v12) {};

\path (15:1.931cm) node[draw=green,shape=circle] (v13) {};
\path (30:1.732cm) node[draw=green,shape=circle] (v14) {};
\draw
(v0) -- (v1)
(v1) -- (v2)
(v2) -- (v0)
(v0) -- (v3)
(v3) -- (v4)
(v4) -- (v0)
(v0) -- (v5)
(v5) -- (v6)
(v6) -- (v0)
(v0) -- (v7)
(v7) -- (v8)
(v8) -- (v0)

(v5) edge[red] (v11)
(v5) edge[red] (v12)
(v11) edge[red] (v12)

(v1) edge[green] (v13)
(v1) edge[green] (v14)
(v13) edge[green] (v14)
%%%%%%%%%%%%%%%%%%%%
(v13) edge[blue] (v9)
(v13) edge[blue] (v10)
(v9) edge[blue] (v10);
\end{tikzpicture}}

%%%%%%%%%%%%%%%%%%%%%%%%%%%%%%%%%%%%% Graph e %%%%%%%%%%%%%%%%%%%%%%%%%%%%%%%%%%%%%%%%%%%%%%%%%%%%%%%%%%%%%%%%%%%%%%%%%%%%%%%%%%%%%
\put(25,95){\begin{tikzpicture}[scale=1.1]
%\tikzstyle{every node}=[draw,shape=circle];
\path (0:0cm) node[draw,shape=circle] (v0) {};
\path (0:1cm) node[draw,shape=circle] (v1) {};
\path (40:1cm) node[draw,shape=circle] (v2) {};
\path (92:1cm) node[draw,shape=circle] (v3) {};
\path (132:1cm) node[draw,shape=circle] (v4) {};
\path (174:1cm) node[draw,shape=circle] (v5) {};
\path (214:1cm) node[draw,shape=circle] (v6) {};
\path (266:1cm) node[draw,shape=circle] (v7) {};
\path (306:1cm) node[draw,shape=circle] (v8) {};
%%%%%%%%%%%%%%%%%%%%%%%%%%%%%%%%%%%%%%%%%%%%%%%%%%%%%%%%%%%%%%%%%%%%%%%%%%%%%%%%
\path (45:1.931cm) node[draw=blue,shape=circle] (v9) {};
\path (60:1.931cm) node[draw=blue,shape=circle] (v10){};

\path (345:1.931cm) node[draw=red,shape=circle] (v11) {};
\path (0:2cm) node[draw=red,shape=circle] (v12) {};

\path (15:1.931cm) node[draw=green,shape=circle] (v13) {};
\path (30:1.732cm) node[draw=green,shape=circle] (v14) {};
\draw
(v0) -- (v1)
(v1) -- (v2)
(v2) -- (v0)
(v0) -- (v3)
(v3) -- (v4)
(v4) -- (v0)
(v0) -- (v5)
(v5) -- (v6)
(v6) -- (v0)
(v0) -- (v7)
(v7) -- (v8)
(v8) -- (v0)

(v1) edge[red] (v11)
(v1) edge[red] (v12)
(v11) edge[red] (v12)

(v1) edge[green] (v13)
(v1) edge[green] (v14)
(v13) edge[green] (v14)
%%%%%%%%%%%%%%%%%%%%
(v9) edge[blue] (v2)
(v10) edge[blue] (v2)
(v9) edge[blue] (v10);
\end{tikzpicture}}
%%%%%%%%%%%%%%%%%%%%%%%%%%%%%%%%%%%%% Graph f %%%%%%%%%%%%%%%%%%%%%%%%%%%%%%%%%%%%%%%%%%%%%%%%%%%%%%%%%%%%%%%%%%%%%%%%%%%%%%%%%%%%%
\put(220,95){\begin{tikzpicture}[scale=1.1]
%\tikzstyle{every node}=[draw,shape=circle];
\path (0:0cm) node[draw,shape=circle] (v0) {};
\path (0:1cm) node[draw,shape=circle] (v1) {};
\path (40:1cm) node[draw,shape=circle] (v2) {};
\path (92:1cm) node[draw,shape=circle] (v3) {};
\path (132:1cm) node[draw,shape=circle] (v4) {};
\path (174:1cm) node[draw,shape=circle] (v5) {};
\path (214:1cm) node[draw,shape=circle] (v6) {};
\path (266:1cm) node[draw,shape=circle] (v7) {};
\path (306:1cm) node[draw,shape=circle] (v8) {};
%%%%%%%%%%%%%%%%%%%%%%%%%%%%%%%%%%%%%%%%%%%%%%%%%%%%%%%%%%%%%%%%%%%%%%%%%%%%%%%%
%\path (30:2.631cm) node[draw=blue,shape=circle] (v9) {};
%\path (40:2.332cm) node[draw=blue,shape=circle] (v10) {};
\path (30:1.931cm) node[draw=blue,shape=circle] (v9) {};
\path (45:1.732cm) node[draw=blue,shape=circle] (v10) {};

\path (-15:1.931cm) node[draw=red,shape=circle] (v11) {};
\path (-30:2cm) node[draw=red,shape=circle] (v12) {};

\path (0:2.131cm) node[draw=green,shape=circle] (v13) {};
\path (15:1.932cm) node[draw=green,shape=circle] (v14) {};
\draw
(v0) -- (v1)
(v1) -- (v2)
(v2) -- (v0)
(v0) -- (v3)
(v3) -- (v4)
(v4) -- (v0)
(v0) -- (v5)
(v5) -- (v6)
(v6) -- (v0)
(v0) -- (v7)
(v7) -- (v8)
(v8) -- (v0)

(v1) edge[red] (v11)
(v1) edge[red] (v12)
(v11) edge[red] (v12)

(v1) edge[green] (v13)
(v1) edge[green] (v14)
(v13) edge[green] (v14)
%%%%%%%%%%%%%%%%%%%%
(v1) edge[blue] (v9)
(v1) edge[blue] (v10)
(v9) edge[blue] (v10);
\end{tikzpicture}}
%%%%%%%%%%%%%%%%%%%%%%%%%%%%%%%%%%%%%%%%%%%%%%%%%%%%%%%%%%%%%%%%%%%%%%%%%%%%%%%%%%%%%%%%%%%%%%%%%%%%%%%%%%%%%%%%%%%%%%%%%%

%%%%%%%%%%%%%%%%%%%%%%%%%%%%%%%%%%%%%%%%%%%%%%%%%%%%%%%%%%%%%%%%%%%%%%%%%%%%%%%%%%%%%%%%%%%%%%%%%%%%%%%%%%%%%%%%%%%%%%%%%%%%%%%%%%%%%

\end{picture}
%\caption{All possible Dutch windmills of order $7$ with two pendant triangles.}
\label{twopendantfig}
\end{center}
\end{figure}

%%%%%%%%%%%%%%%%%%%%%%%%%%%%%%%%%%%%%%%%%%%
%%%%%%%%%%%%%%%%%%%%%%%%%%%%%%%%%%%%%%%%%%%%%%%
%%%%%%%%%%%%%%%%%%%%%%%%%%%%%%%%%%%%%%%%%%%%%%%%%%%%%%%%
%%%%%%%%%%%%%%%%%%%%%%%%%%%%%%%%%%%%%%%%%%%%%%%%%%%%%%%%%%%
\begin{figure}[H]
\begin{center}
\begin{picture}(510,300)
\linethickness{0.9pt}
%\put(25,630){$(a)$}     %
%\put(220,630){$(b)$}    %

%\put(25,510){$(c)$}    %
%\put(220,510){$(d)$}  %
  %
\put(25,220){$(g)$}  %  300,860
\put(220,220){$(h)$} % 25,860

\put(25,120){$(i)$}    %
\put(220,120){$(j)$}  %

\put(122,5){$(k)$}  %

%%%%%%%%%%%%%%%%%%%%%%%%%%%%%%%%%%%%%%%

%%%%%%%%%%%%%%%%%%%%%%%%%%%%%%%%%%%%% Graph g %%%%%%%%%%%%%%%%%%%%%%%%%%%%%%%%%%%%%%%%%%%%%%%%%%%%%%%%%%%%%%%%%%%%%%%%%%%%%%%%%%%%%

\put(25,230){\begin{tikzpicture}[scale=1.1]
%\tikzstyle{every node}=[draw,shape=circle];
\path (0:0cm) node[draw,shape=circle] (v0) {};
\path (0:1cm) node[draw,shape=circle] (v1) {};
\path (40:1cm) node[draw,shape=circle] (v2) {};
\path (92:1cm) node[draw,shape=circle] (v3) {};
\path (132:1cm) node[draw,shape=circle] (v4) {};
\path (174:1cm) node[draw,shape=circle] (v5) {};
\path (214:1cm) node[draw,shape=circle] (v6) {};
\path (266:1cm) node[draw,shape=circle] (v7) {};
\path (306:1cm) node[draw,shape=circle] (v8) {};
%%%%%%%%%%%%%%%%%%%%%%%%%%%%%%%%%%%%%%%%%%%%%%%%%%%%%%%%%%%%%%%%%%%%%%%%%%%%%%%%
\path (45:1.931cm) node[draw=blue,shape=circle] (v9) {};
\path (60:1.931cm) node[draw=blue,shape=circle] (v10){};

\path (15:2.931cm) node[draw=red,shape=circle] (v11) {};
\path (30:2.732cm) node[draw=red,shape=circle] (v12) {};

\path (15:1.931cm) node[draw=green,shape=circle] (v13) {};
\path (30:1.732cm) node[draw=green,shape=circle] (v14) {};
\draw
(v0) -- (v1)
(v1) -- (v2)
(v2) -- (v0)
(v0) -- (v3)
(v3) -- (v4)
(v4) -- (v0)
(v0) -- (v5)
(v5) -- (v6)
(v6) -- (v0)
(v0) -- (v7)
(v7) -- (v8)
(v8) -- (v0)

(v13) edge[red] (v11)
(v13) edge[red] (v12)
(v11) edge[red] (v12)

(v1) edge[green] (v13)
(v1) edge[green] (v14)
(v13) edge[green] (v14)
%%%%%%%%%%%%%%%%%%%%
(v9) edge[blue] (v2)
(v10) edge[blue] (v2)
(v9) edge[blue] (v10);
\end{tikzpicture}}

%%%%%%%%%%%%%%%%%%%%%%%%%%%%%%%%%%%%% Graph h %%%%%%%%%%%%%%%%%%%%%%%%%%%%%%%%%%%%%%%%%%%%%%%%%%%%%%%%%%%%%%%%%%%%%%%%%%%%%%%%%%%%%

\put(220,230){\begin{tikzpicture}[scale=1.1]
%\tikzstyle{every node}=[draw,shape=circle];
\path (0:0cm) node[draw,shape=circle] (v0) {};
\path (0:1cm) node[draw,shape=circle] (v1) {};
\path (40:1cm) node[draw,shape=circle] (v2) {};
\path (92:1cm) node[draw,shape=circle] (v3) {};
\path (132:1cm) node[draw,shape=circle] (v4) {};
\path (174:1cm) node[draw,shape=circle] (v5) {};
\path (214:1cm) node[draw,shape=circle] (v6) {};
\path (266:1cm) node[draw,shape=circle] (v7) {};
\path (306:1cm) node[draw,shape=circle] (v8) {};
%%%%%%%%%%%%%%%%%%%%%%%%%%%%%%%%%%%%%%%%%%%%%%%%%%%%%%%%%%%%%%%%%%%%%%%%%%%%%%%%
\path (30:2.631cm) node[draw=blue,shape=circle] (v9) {};
\path (40:2.332cm) node[draw=blue,shape=circle] (v10) {};

\path (0:2.931cm) node[draw=red,shape=circle] (v11) {};
\path (15:2.732cm) node[draw=red,shape=circle] (v12) {};

\path (15:1.931cm) node[draw=green,shape=circle] (v13) {};
\path (30:1.732cm) node[draw=green,shape=circle] (v14) {};
\draw
(v0) -- (v1)
(v1) -- (v2)
(v2) -- (v0)
(v0) -- (v3)
(v3) -- (v4)
(v4) -- (v0)
(v0) -- (v5)
(v5) -- (v6)
(v6) -- (v0)
(v0) -- (v7)
(v7) -- (v8)
(v8) -- (v0)

(v13) edge[red] (v11)
(v13) edge[red] (v12)
(v11) edge[red] (v12)

(v1) edge[green] (v13)
(v1) edge[green] (v14)
(v13) edge[green] (v14)
%%%%%%%%%%%%%%%%%%%%
(v13) edge[blue] (v9)
(v13) edge[blue] (v10)
(v9) edge[blue] (v10);
\end{tikzpicture}}

%%%%%%%%%%%%%%%%%%%%%%%%%%%%%%%%%%%%% Graph i %%%%%%%%%%%%%%%%%%%%%%%%%%%%%%%%%%%%%%%%%%%%%%%%%%%%%%%%%%%%%%%%%%%%%%%%%%%%%%%%%%%%%

\put(25,130){\begin{tikzpicture}[scale=1.1]
%\tikzstyle{every node}=[draw,shape=circle];
\path (0:0cm) node[draw,shape=circle] (v0) {};
\path (0:1cm) node[draw,shape=circle] (v1) {};
\path (40:1cm) node[draw,shape=circle] (v2) {};
\path (92:1cm) node[draw,shape=circle] (v3) {};
\path (132:1cm) node[draw,shape=circle] (v4) {};
\path (174:1cm) node[draw,shape=circle] (v5) {};
\path (214:1cm) node[draw,shape=circle] (v6) {};
\path (266:1cm) node[draw,shape=circle] (v7) {};
\path (306:1cm) node[draw,shape=circle] (v8) {};
%%%%%%%%%%%%%%%%%%%%%%%%%%%%%%%%%%%%%%%%%%%%%%%%%%%%%%%%%%%%%%%%%%%%%%%%%%%%%%%%
\path (30:2.631cm) node[draw=blue,shape=circle] (v9) {};
\path (40:2.332cm) node[draw=blue,shape=circle] (v10) {};

\path (345:1.931cm) node[draw=red,shape=circle] (v11) {};
\path (0:2cm) node[draw=red,shape=circle] (v12) {};

\path (15:1.931cm) node[draw=green,shape=circle] (v13) {};
\path (30:1.732cm) node[draw=green,shape=circle] (v14) {};
\draw
(v0) -- (v1)
(v1) -- (v2)
(v2) -- (v0)
(v0) -- (v3)
(v3) -- (v4)
(v4) -- (v0)
(v0) -- (v5)
(v5) -- (v6)
(v6) -- (v0)
(v0) -- (v7)
(v7) -- (v8)
(v8) -- (v0)

(v1) edge[red] (v11)
(v1) edge[red] (v12)
(v11) edge[red] (v12)

(v1) edge[green] (v13)
(v1) edge[green] (v14)
(v13) edge[green] (v14)
%%%%%%%%%%%%%%%%%%%%
(v13) edge[blue] (v9)
(v13) edge[blue] (v10)
(v9) edge[blue] (v10);
\end{tikzpicture}}

%%%%%%%%%%%%%%%%%%%%%%%%%%%%%%%%%%%%% Graph j %%%%%%%%%%%%%%%%%%%%%%%%%%%%%%%%%%%%%%%%%%%%%%%%%%%%%%%%%%%%%%%%%%%%%%%%%%%%%%%%%%%%
\put(220,130){\begin{tikzpicture}[scale=1.1]
%\tikzstyle{every node}=[draw,shape=circle];
\path (0:0cm) node[draw,shape=circle] (v0) {};
\path (0:1cm) node[draw,shape=circle] (v1) {};
\path (40:1cm) node[draw,shape=circle] (v2) {};
\path (92:1cm) node[draw,shape=circle] (v3) {};
\path (132:1cm) node[draw,shape=circle] (v4) {};
\path (174:1cm) node[draw,shape=circle] (v5) {};
\path (214:1cm) node[draw,shape=circle] (v6) {};
\path (266:1cm) node[draw,shape=circle] (v7) {};
\path (306:1cm) node[draw,shape=circle] (v8) {};
%%%%%%%%%%%%%%%%%%%%%%%%%%%%%%%%%%%%%%%%%%%%%%%%%%%%%%%%%%%%%%%%%%%%%%%%%%%%%%%%
\path (30:2.631cm) node[draw=blue,shape=circle] (v9) {};
\path (40:2.332cm) node[draw=blue,shape=circle] (v10) {};

\path (20:3.331cm) node[draw=red,shape=circle] (v11) {};
\path (30:3.632cm) node[draw=red,shape=circle] (v12) {};

\path (15:1.931cm) node[draw=green,shape=circle] (v13) {};
\path (30:1.732cm) node[draw=green,shape=circle] (v14) {};
\draw
(v0) -- (v1)
(v1) -- (v2)
(v2) -- (v0)
(v0) -- (v3)
(v3) -- (v4)
(v4) -- (v0)
(v0) -- (v5)
(v5) -- (v6)
(v6) -- (v0)
(v0) -- (v7)
(v7) -- (v8)
(v8) -- (v0)

(v9) edge[red] (v11)
(v9) edge[red] (v12)
(v11) edge[red] (v12)

(v1) edge[green] (v13)
(v1) edge[green] (v14)
(v13) edge[green] (v14)
%%%%%%%%%%%%%%%%%%%%
(v13) edge[blue] (v9)
(v13) edge[blue] (v10)
(v9) edge[blue] (v10);
\end{tikzpicture}}
%%%%%%%%%%%%%%%%%%%%%%%%%%%%%%%%%%%%% Graph k %%%%%%%%%%%%%%%%%%%%%%%%%%%%%%%%%%%%%%%%%%%%%%%%%%%%%%%%%%%%%%%%%%%%%%%%%%%%%%%%%%%%%
\put(122,15){\begin{tikzpicture}[scale=1.1]
%\tikzstyle{every node}=[draw,shape=circle];
\path (0:0cm) node[draw,shape=circle] (v0) {};
\path (0:1cm) node[draw,shape=circle] (v1) {};
\path (40:1cm) node[draw,shape=circle] (v2) {};
\path (92:1cm) node[draw,shape=circle] (v3) {};
\path (132:1cm) node[draw,shape=circle] (v4) {};
\path (174:1cm) node[draw,shape=circle] (v5) {};
\path (214:1cm) node[draw,shape=circle] (v6) {};
\path (266:1cm) node[draw,shape=circle] (v7) {};
\path (306:1cm) node[draw,shape=circle] (v8) {};
%%%%%%%%%%%%%%%%%%%%%%%%%%%%%%%%%%%%%%%%%%%%%%%%%%%%%%%%%%%%%%%%%%%%%%%%%%%%%%%%
\path (3:2.631cm) node[draw=blue,shape=circle] (v9) {};
\path (-12:2.332cm) node[draw=blue,shape=circle] (v10) {};
\path (27:2.631cm) node[draw=red,shape=circle] (v11) {};
\path (40:2.732cm) node[draw=red,shape=circle] (v12) {};
\path (15:1.931cm) node[draw=green,shape=circle] (v13) {};
\path (30:1.732cm) node[draw=green,shape=circle] (v14) {};
\draw
(v0) -- (v1)
(v1) -- (v2)
(v2) -- (v0)
(v0) -- (v3)
(v3) -- (v4)
(v4) -- (v0)
(v0) -- (v5)
(v5) -- (v6)
(v6) -- (v0)
(v0) -- (v7)
(v7) -- (v8)
(v8) -- (v0)

(v14) edge[red] (v11)
(v14) edge[red] (v12)
(v11) edge[red] (v12)

(v1) edge[green] (v13)
(v1) edge[green] (v14)
(v13) edge[green] (v14)
%%%%%%%%%%%%%%%%%%%%
(v13) edge[blue] (v9)
(v13) edge[blue] (v10)
(v9) edge[blue] (v10);
\end{tikzpicture}}
\end{picture}
\end{center}
\caption{All types of Dutch windmills of order $7$ with three pendant triangles.}\label{threependantfig}
\end{figure}
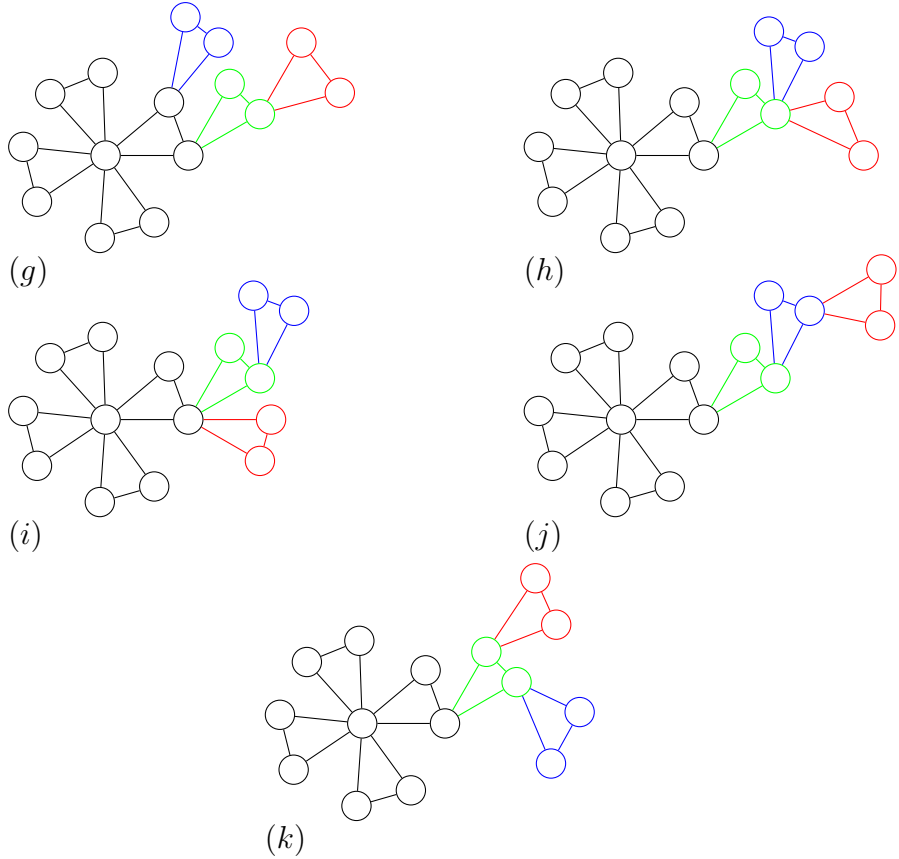
%\rotatebox[origin=c]{90}{Figure~\thefigure: My rotated caption}
%\caption[My short caption]{\lipsum[2]}
%%%%%%%%%%%%%%%%%%%%%%%%%%%%%%%%%%%%%%%%%%%%%%%%%%%%%%%%%%%%%%%%%%%%%%%%%%%%%%%%%%%%%%%%%%%%%%%%%%%%%%%%%%%%%%%%%%%%%%%%%%%%%%%%%%%%%%%%%%%%%%%%%%%%%%%%%%%%%%%%%%%%%%%%%%%%%%%%%%%%%%%%%%%%%%%%%%%%%%%%%%%%%%%%%%%%%%%%%%%%%%%%%%%%%%%%%%%%%%%%%%%%%%%%%%%%%%%%%%%%%%%%%%%%%%%%%%%%%%%%%%%%%%%%%%%%%%%%%%%%%%%%
\begin{theorem}\label{dutchthrm1}\cite{Bermond} Let $G$ be a Dutch windmill with $n$ blocks. If there exists a Skolem sequence of order $n$, then $G$ is graceful.\end{theorem}
\begin{pf}Let $G$ be a Dutch windmill with $n$ blocks. Let $S_{n}$ be a Skolem sequence of order $n$ of the form $\left(a_{i},b_{i}\right)$ with $a_{i} < b_{i},$ for $i=1,2,\ldots,n.$ These pairs give $n$ base blocks which are $\{0,a_{i}+n,b_{i}+n\}$, $1\leq i\leq n.$
\begin{enumerate}
	\item We can obtain three types of differences from $\{0,a_{i}+n,b_{i}+n\}$ as follows:
	\begin{enumerate}
		\item $A=\left\{(b_{i}+n)-(a_{i}+n)\right\}=\left\{1,2,\ldots,n\right\}$;
		\item $B=\left\{(a_{i}+n)-0\right\}$;
		\item $C=\left\{(b_{i}+n)-0\right\}$.
	\end{enumerate}
Then, $A\cup B \cup C=\left\{1,2,\ldots,3n\right\}$.
	\item The base blocks $\{0,a_i+n,b_i+n\}_{i=1}^n$ give the vertex labels $\{0$, $1$, $2$, $\ldots$, $3n\},$ where $0$ is repeated $n$ times, but it is a common vertex. The base blocks formed by Skolem sequences present a graceful labelling.\end{enumerate}\end{pf}
%ttttttttttttttttttttttttttttttttttttttttttttttttttttttttttttttttttttttttttttttttttttttttttttttttttttttttttttttttttttttttttttttttttttttttttttttttttttttttttttttttttttttttttttttttttttttttttttttttttttttttttttttttttttttttttttttttttttttttttttttttttttttttttttttttttttttttttttttttttttttttttttttttttttttttttttttt
\begin{theorem}\label{dutchthrm2}\cite{Bermond} Let $G$ be a Dutch windmill with $n$ blocks. If there exists a hooked Skolem sequence of order $n$, then $G$ is near graceful.\end{theorem}
\begin{pf}Let $G$ be a Dutch windmill with $n$ blocks. Let $hS_{n}$ be a hooked Skolem sequence of order $n$ of the form $\left(a_{i},b_{i}\right)$ with $a_{i} < b_{i},$ for $i=1,2,\ldots,n.$ These pairs give $n$ base blocks which are $\{0,a_{i}+n,b_{i}+n\}$, $1\leq i\leq n.$
\begin{enumerate}
	\item We can obtain three types of differences from $\{0,a_{i}+n,b_{i}+n\}$ as follows:
	\begin{enumerate}
		\item $A=\left\{(b_{i}+n)-(a_{i}+n)\right\}=\left\{1,2,\ldots,n\right\}$;
		\item $B=\left\{(a_{i}+n)-0\right\}$;
		\item $C=\left\{(b_{i}+n)-0\right\}$.
	\end{enumerate}
Then, $A\cup B \cup C=\left\{1,2,\ldots,3n-1,3n+1\right\}$.
	\item The base blocks $\{0,a_i+n,b_i+n\}_{i=1}^n$ give the vertex labels $\{0, 1, 2$, $\ldots$, $3n-1$, $3n+1\}$, where $0$ is repeated $n$ times, but it is a common vertex. The base blocks formed by a hooked Skolem sequences present a near graceful labelling.\end{enumerate} \end{pf}
%\begin{theorem}\label{dutchthrm}Take the elements of the base blocks that constructed from a (hooked) Skolem sequence $(0,a_i+n,b_i+n)$ or $(0,i,b_i+n)$ formed by (hooked) Skolem sequences these elements give the vertex labels of the (near) gracefully labelled Dutch windmill.\end{theorem}
%\par Furthermore, take the differences between the elements of the base blocks that constructed from a (hooked) Skolem sequence $(0,a_i+n,b_i+n)$ or $(0,i,b_i+n)$ formed by (hooked) Skolem sequences these differences give the edges labels of the (near) gracefully labelled Dutch windmill.\end{remark}\par
In the (near) graceful labelling by (hooked) Skolem sequences, we can use the base blocks of the form $\{0,a_i+n,b_i+n\}_{i=1}^n$ or $\{0,i,b_i+n\}_{i=1}^n$. They give two different vertex labels, but the same edge labels.
%%%%%%%%%%%%%%%%%%%%%%%%%%%%%%%%%%%%%%%%%%%%%%%%%%%%%%%%%%%%%%%%%%%%%%%%%%%%%%%%%%%%%%%%%%%%%%%%%%%%%%%%%%%%%%%%%%%%%%%%%%%%
%Here we will present how we can graceful labelling a Dutch windmill containing four triangles. Start with a Skolem sequence of order four and form the triples of the form $\{0,a_i+n,b_i+n\}$ or $\{0,i,b_i+n\}$. These triples gracefully label the Dutch windmill order 4. By the constructions in Example~\ref{baseblocks} we can see a three different gracefully labelled of the Dutch windmill of order 4.\par
For instance, taking the base blocks from the example in Section~\ref{skolemintro} will give us the gracefully labelled Dutch windmill of order $4$ shown in Figures~\ref{different3}$(a),(b),$ and $(c).$ The gracefully labelled Dutch windmills in Figure \ref{different3}$(a)$ are formed by base blocks of the form $\left\{0,a_i+n,b_i+n\right\}$; the base blocks of the form $\left\{0,i,b_i+n\right\}$ in Figure \ref{different3}$(b)$; and finally by a mixed set of base blocks that come from both forms in Figure \ref{different3}$(c)$. We will use this technique in Section 2 when creating graceful labellings for Type $(k)$ Dutch windmills.
\\
\\
\begin{figure}[h!]
\begin{center}
\begin{picture}(260,235)
\linethickness{0.9pt}
\put(0,148){$(a)$}     %
\put(155,148){$(b)$}    %
\put(78,0){$(c)$}    %
%\begin{tikzpicture}[scale=1.3]
\put(0,150){\begin{tikzpicture}[scale=1.2]
%\tikzstyle{every node}=[draw,shape=circle];
\path (0:0cm) node[draw,shape=circle] (v0) {$0$};
\path (20:1.5cm) node[draw,shape=circle] (v1) {\tiny $5$};
\path (70:1.5cm) node[draw,shape=circle] (v2) {\tiny $6$};
\path (110:1.5cm) node[draw,shape=circle] (v3) {\tiny $8$};
\path (160:1.5cm) node[draw,shape=circle] (v4) {\tiny $10$};
\path (200:1.5cm) node[draw,shape=circle] (v5) {\tiny $9$};
\path (250:1.5cm) node[draw,shape=circle] (v6) {\tiny $12$};
\path (290:1.5cm) node[draw,shape=circle] (v7) {\tiny $7$};
\path (340:1.5cm) node[draw,shape=circle] (v8) {\tiny $11$};
\draw
(v0) -- (v1) node[pos=.6,below] {\tiny $5$}
(v0) -- (v2) node[pos=.6,right] {\tiny $6$}
(v0) -- (v3) node[pos=.5,left] {\tiny $8$}
(v0) -- (v4) node[pos=.6,below] {\tiny $10$}
(v0) -- (v5) node[pos=.5,below] {\tiny $9$}
(v0) -- (v6) node[pos=.7,right] {\tiny $12$}
(v0) -- (v7) node[pos=.5,right] {\tiny $7$}
(v0) -- (v8) node[pos=.5,below] {\tiny $11$}
(v2) -- (v1) node[pos=.4,below] {\tiny $1$}
(v3) -- (v4) node[pos=.4,below] {\tiny $2$}
(v5) -- (v6) node[pos=.3,below] {\tiny $3$}
(v7) -- (v8) node[pos=.6,below] {\tiny $4$};
\end{tikzpicture}}
%\end{center}
%\begin{remark}[Note to myself]\textbf{I think here we need to define a triangular snake with a pendant}\end{remark}

\put(155,150){\begin{tikzpicture}[scale=1.2]
%\tikzstyle{every node}=[draw,shape=circle];
\path (0:0cm) node[draw,shape=circle] (v0) {\tiny $0$};
\path (20:1.5cm) node[draw,shape=circle] (v1) {\tiny $1$};
\path (70:1.5cm) node[draw,shape=circle] (v2) {\tiny $6$};
\path (110:1.5cm) node[draw,shape=circle] (v3) {\tiny $2$};
\path (160:1.5cm) node[draw,shape=circle] (v4) {\tiny $10$};
\path (200:1.5cm) node[draw,shape=circle] (v5) {\tiny $3$};
\path (250:1.5cm) node[draw,shape=circle] (v6) {\tiny $12$};
\path (290:1.5cm) node[draw,shape=circle] (v7) {\tiny $4$};
\path (340:1.5cm) node[draw,shape=circle] (v8) {\tiny $11$};
\draw
(v0) -- (v1) node[pos=.6,below] {\tiny $1$}
(v0) -- (v2) node[pos=.6,right] {\tiny $6$}
(v0) -- (v3) node[pos=.5,left] {\tiny $2$}
(v0) -- (v4) node[pos=.6,below] {\tiny $10$}
(v0) -- (v5) node[pos=.5,below] {\tiny $3$}
(v0) -- (v6) node[pos=.7,right] {\tiny $12$}
(v0) -- (v7) node[pos=.5,right] {\tiny $4$}
(v0) -- (v8) node[pos=.5,below] {\tiny $11$}
(v2) -- (v1) node[pos=.4,below] {\tiny $5$}
(v3) -- (v4) node[pos=.4,below] {\tiny $8$}
(v5) -- (v6) node[pos=.3,below] {\tiny $9$}
(v7) -- (v8) node[pos=.6,below] {\tiny $7$};
\end{tikzpicture}}
%\end{center}

\put(78,5){\begin{tikzpicture}[scale=1.2]
%\tikzstyle{every node}=[draw,shape=circle];
\path (0:0cm) node[draw,shape=circle] (v0) {\tiny $0$};
\path (20:1.5cm) node[draw,shape=circle] (v1) {\tiny $1$};
\path (70:1.5cm) node[draw,shape=circle] (v2) {\tiny $6$};
\path (110:1.5cm) node[draw,shape=circle] (v3) {\tiny $2$};
\path (160:1.5cm) node[draw,shape=circle] (v4) {\tiny $10$};
\path (200:1.5cm) node[draw,shape=circle] (v5) {\tiny $9$};
\path (250:1.5cm) node[draw,shape=circle] (v6) {\tiny $12$};
\path (290:1.5cm) node[draw,shape=circle] (v7) {\tiny $7$};
\path (340:1.5cm) node[draw,shape=circle] (v8) {\tiny $11$};
\draw
(v0) -- (v1) node[pos=.6,below] {\tiny $1$}
(v0) -- (v2) node[pos=.6,right] {\tiny $6$}
(v0) -- (v3) node[pos=.5,left] {\tiny $2$}
(v0) -- (v4) node[pos=.6,below] {\tiny $10$}
(v0) -- (v5) node[pos=.5,below] {\tiny $9$}
(v0) -- (v6) node[pos=.7,right] {\tiny $12$}
(v0) -- (v7) node[pos=.5,right] {\tiny $7$}
(v0) -- (v8) node[pos=.5,below] {\tiny $11$}
(v2) -- (v1) node[pos=.4,below] {\tiny $5$}
(v3) -- (v4) node[pos=.4,below] {\tiny $8$}
(v5) -- (v6) node[pos=.3,below] {\tiny $3$}
(v7) -- (v8) node[pos=.6,below] {\tiny $4$};
\end{tikzpicture}}
\end{picture}
\caption{Three different gracefully labelled Dutch windmills of order 4.}\label{different3}
\end{center}
\end{figure}
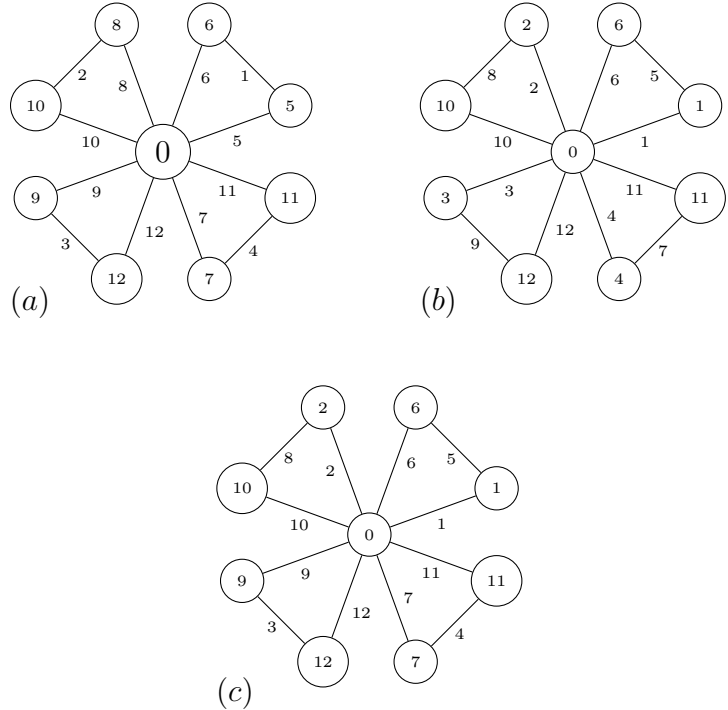
\par Since $x-y=(x+c)-(y+c)$ for any $x,y$, the following lemma is straightforward.
\begin{lemma}\label{lemma2}If we add a constant $c$ to each element of any triple $\left\{x,y,z\right\},$ which results in $\left\{x+c,y+c,z+c\right\}$, then the differences between the elements of $\left\{x+c,y+c,z+c\right\}$ will be the same.\end{lemma}
%If we have a pivot triple of the form $\left(0,a_{i}+n,b_{i}+n\right)$ and replace it with\newline $\left(i,a_{i}+i+n,b_{i}+i+n\right)$, with $1 \leq i \leq n$ then we call this technique \textit{pivoting}.
In the special case of having a triple of the form $\left(0,a_{i}+n,b_{i}+n\right)$, where $a_i$ is a pivot, and replacing that triple with the new triple $\left(i,a_{i}+i+n,b_{i}+i+n\right)$, we say that we are \textit{pivoting}.

Consider the triples as given in Figure~\ref{different3} (a), then add $2$ to each element of the triple $\{0,8,10\}$ to obtain $\{2,10,12\}$. This gives the blocks $\{0,5,6\}$, $\{2,10,12\}$, $\{0,9,12\}$, and $\{0,7,11\}.$ Consequently, we have obtained a graceful labelling for a Dutch windmill of order three with one pendant triangle.

Here we will present an example of how to label Figure~\ref{threependantfig}(b) by using a hooked Skolem sequence.
{\label{threepivotsexample}}Let $hS_{7}=(7$, 4, 6, 3, 5, 4, 3, 7, 6, 5, 1, 1, 2, 0, $2)$ be a hooked Skolem sequence of order $7.$ This yields the pairs $\{(11,12)$, $(13,15)$, $(4,7)$, $(2,6)$, $(5,10)$, $(3,9)$, $(1,8)\}$. This sequence gives the base blocks of the form $\left\{0,a_{i}+n, b_{i}+n\right\}$ as follows: $A=\{(0,18,19)$, $(0,20,22)$, $(0,11,14)$, $(0,9,13)$, $(0,12,17)$, $(0,10,16)$, $(0,8,15)\}$. The above sequence $hS_{7}$ has six pivots, but we consider only three in particular, which are $1,4,$ and $7.$ We obtain a new set of base blocks by pivoting as follows: $A'=\{(1,19,20)$, $(0,20,22)$, $(0,11,14)$, $(4,13,17)$, $(0,12,17)$, $(0,10,16)$, $(7,15,22)\}.$

We will use the set $A'$ to label our graph as follows: label the central vertex with $0$. Then label all the vanes by these base blocks containing $0.$ Finally, label the pendant triangles with the triples corresponding to the pivots. The blocks  $\left\{1,19,20\right\}$ and $\left\{7,15,22\right\}$ each intersect $\left\{0,20,22\right\}$ at a single distinct element, namely $20$ for the first block and $22$ for the second block. The block $\left\{4,13,17\right\}$ intersects $\left\{0,12,17\right\}$ at a single element, namely $17.$

To label Dutch windmills with $n$ blocks that have three pendant triangles, we will use Skolem sequences of order $n$ with at least three pivots. Since our Skolem sequence $S$ has three distinct pivots $r,s,t,\  1 \leq r,s,t \leq n$ we will obtain a new set of base blocks $A'$. Then $A'$ contains the new triples $\{r,a_{r}+n+r,b_{r}+n+r\}, \{s,a_{s}+n+s, b_{s}+n+s\}, \{t,a_{t}+n+t,b_{t}+n+t\}.$ Exactly one of $r,a_{r}+n+r,$ or $b_{r}+n+r$ must appear in one of the other triples and this is the label where we attach the pendant triangle; the elements of the triples in the set $A'$ must be labelled from  $\left\{0,1,2,\ldots,3n\right\}$, and these values cannot be more than $3n$, and a similar argument holds for $s$ and $t$. Then, by Theorem~\ref{dutchthrm1}, the Skolem sequence gives the vertex labels of the gracefully labelled Dutch windmill. Likewise, the same idea works for hooked Skolem sequences, and, by Theorem~\ref{dutchthrm2}, the hooked Skolem sequence gives the vertex labels of the near gracefully labelled Dutch windmill.
%%%%%%%%%%%%%%%%%%%%%%%%%%%%%%%%%%%%%%%%%%%%%%%%%%%%%%%%%%%%%%%%%%%%%%%%%%%%%%%%%%%%%%%%%%%%%%%%%%%%%%%%%%%%%%%%%%%%%%%%%%%%%%%%%%%%%%%%%%%%%%%%%%%%%%%%

In the following lemma, $\lambda$ represents a letter between $a$ and $k.$
\begin{lemma}\label{bound}If a (hooked) Skolem sequence of order $n$ exists that gives a (near) graceful labelling of a Type $(\lambda)$ Dutch windmill with $n$ blocks then for $m \geq 3n+1,$ there exists a (near) graceful labelling of a Type $(\lambda)$ Dutch windmill with $m$ blocks. \end{lemma}
\begin{pf}Let $G$ be a Type $(\lambda)$ Dutch windmill with $n$ blocks (near) gracefully labelled by $S_{n},$ a (hooked) Skolem sequence of order $n$. By Theorem~\ref{langfordexistence}, a (hooked) Langford sequence of order $m$ exists with defect $d=n+1.$ Then by Lemma~\ref{skla}, we obtain a new (hooked) Skolem sequence of order $m,$ $S_{m}$ where $m \geq 3n + 1=N(n)$. Then, since $S_n$ had the pivot structure needed to (near) gracefully label a Type $(\lambda)$ Dutch windmill with $n$ blocks, $S_m$ will yield the structure needed to (near) gracefully label a Type $(\lambda)$ Dutch windmill with $m$ blocks.
\end{pf}
% The sequence $S_{m}$ gives a graceful labelling of a dutch windmill with $m$ blocks.Then by the structure of the sequence $S_{n},$ we use to label $G$ we get a (near) graceful label of Type $(\lambda)$ dutch windmill with $m$ blocks.
%\textbf{we need to talk about how we can label and where we can attach the pivot.}
%\end{proof}
%\pagebreak
%%%%%%%%%%%%%%%%%%%%%%%%%%%%%%%%%%%%%%%%%%%%%%%%%%%%%%%%%%%%%%%%%%%%%%%%%%%%%%%%%%%%%%%%%%%%%%%%%%%%%%%%%%%%%%%%%%%%%%%%%%%%%%%%%%%%%%%%%%%%%%%%%%%%%%%%%%%%%%%%%%%%%%%%%%%%%%%%%%%%%%%%%%%%%%%%%%%%%%%%%%%%%%%%%%%%%%%%%%%%%%%%%%%%%%%%%%%%%%%%%%%%%%%%%%%%%%%%%%%%%%%%%%%%%%%%%%%%%%%%%%%%%%%%%%%%%%%%%%%%%%%%
%%%%%%%%%%%%%%%%%%%%%%%%%%%%%% {{  Type 1 }} %%%%%%%%%%%%%%%%%%%%%%%%%%%%%%%%%%%%%%%%%%%%%%%%%%%%%%%%%%%%%%%%%%%%%%%%%%%%%%%%%%%%%%%%%%%%%%%%%%%%%%%%%%%
\subsection{Type $(a)$}
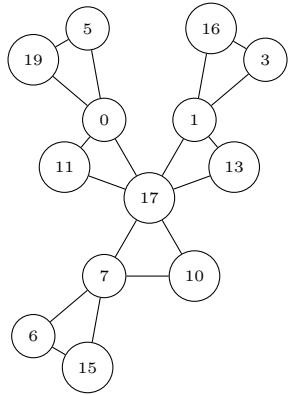
\begin{figure}[H]
\begin{center}
\begin{tikzpicture}[scale=1.2]
%\tikzstyle{every node}=[draw,shape=circle];
\path (0:0cm) node[draw,shape=circle] (v0) {\tiny $17$};
\path (20:1cm) node[draw,shape=circle] (v1) {\tiny $13$};
\path (60:1cm) node[draw,shape=circle] (v2) {\tiny $1$};
\path (120:1cm) node[draw,shape=circle] (v5) {\tiny $0$};
\path (160:1cm) node[draw,shape=circle] (v6) {\tiny $11$};
\path (70:2cm) node[draw,shape=circle] (v11) {\tiny $16$};
\path (50:2cm) node[draw,shape=circle] (v12) {\tiny $3$};
\path (110:2cm) node[draw,shape=circle] (v13) {\tiny $5$};
\path (130:2cm) node[draw,shape=circle] (v14) {\tiny $19$};
\path (90:0.9cm) node[draw=none,shape=circle] (v15) {} ;

\path (-60:1cm) node[draw,shape=circle] (v16) {\tiny $10$};
\path (-120:1cm) node[draw,shape=circle] (v17) {\tiny $7$};

\path (-110:2cm) node[draw,shape=circle] (v18){\tiny $15$};
\path (-130:2cm) node[draw,shape=circle] (v19) {\tiny $6$};
\draw
(v0) -- (v1)
(v1) -- (v2)
(v2) -- (v0)
(v0) -- (v5)
(v5) -- (v6)
(v6) -- (v0)
(v2) -- (v11)
(v2) -- (v12)
(v11) -- (v12)
(v5) -- (v13)
(v5) -- (v14)
(v13) -- (v14)
(v0) -- (v16)
(v0) -- (v17)
(v16) -- (v17)
(v17) -- (v18)
(v17) -- (v19)
(v18) -- (v19);
\end{tikzpicture}
\caption{A near gracefully labelled Type $(a)$ Dutch windmill with $6$ blocks.}\label{taypean6}
\end{center}
\end{figure}

For Type $(a)$, when $n \leq 5$, there are not enough triangles to form a Type $(a)$ Dutch windmill. When $n=6$, see the near graceful labelling in Figure~\ref{taypean6}.
%%%%%%%%%%%%%%%%%%%%%%%%%%%%%%%%%%%%%%%%%%%%%%%%%%%%%%%%%%%%%%%%%%%%%%%%%%%%%%%%%%%%%%%%%%%%%%%%%%%%%%%%%%%%%%%%%%%%%%%%%%%%%%%%%%%%%%%%%%%%%%%%%%%%%%%%%%%%%%%%%%%%%%%%%%%%%%%%%%%%%%%%%%%%%%%%%%%%%%%%%%%%%%%%%%%%%%%%%%%%%%%%%%%%%%%%%%%%%%%%%%%%%%%%%%%%%%%%%%%%%%%%%%%%%%%%%%%%%%%%%%%%%%%%%%%%%%%%%%%%%%%%
\begin{lemma}\label{typea}
Any Type $(a)$ Dutch windmill with at least $6$ blocks is graceful or near graceful.\end{lemma}
\begin{pf}We begin by constructing a (hooked) Skolem sequence which we will use to label a Type $(a)$ Dutch windmill of order 7. Consider the (hooked) Skolem sequence, $hS_{7}=\left(3,4,7,3,2,4,2,5,6,7,1,1,5,0,6\right)$. This sequence gives triples of the form $\left(0,a_{i}+n, b_{i}+n\right)$ as follows: $\left(0,18,19\right)$, $\left(0,12,14\right)$, $\left(0,8,11\right)$, $\left(0,9,13\right)$, $\left(0,15,20\right)$, $\left(0,16,22\right)$, $\left(0,10,17\right)$. The above sequence $hS_{7}$ has three pivots, which are $1,3,$ and $4.$ We convert these triples by pivoting as follows: $\left(1, 19, 20\right)$, $\left(0, 12, 14\right)$, $\left(3, 11, 14\right)$, $\left(4, 13, 17\right)$, $\left(0, 15, 20\right)$, $\left(0, 16, 22\right)$, $\left(0, 10, 17\right)$. These triples near gracefully label the Dutch windmill with $7$ blocks.
%By Step $1,$ we have a hooked Skolem sequence for this type is of order $7, S_{7.}$ Using Lemma~\ref{skla}, we can find a Langford sequence and hooked Langford sequence with defect $d=8.$ We then can obtain a (hooked) Skolem sequence with a considerable order.

We have an order $7$ hooked Skolem sequence. Then by Theorem~\ref{langfordexistence}, a hooked Langford sequence of order $m$ with $d=8$ exists for $m \geq 15$, with $m$ congruent to $1$ and $2$ $({\rm mod}\ 4)$. By Lemma~\ref{skla}, we obtain an associated Skolem sequence of order $n,$ and by Lemma~\ref{bound}, we obtain a gracefully labelled Type $(a)$ Dutch windmill with $n$ blocks for $n \geq 22.$ Likewise, By Theorem~\ref{langfordexistence}, a Langford sequence of order $m$ with $d=8$ exists for $m \geq 15$ with $m$ congruent to $0$ and $3$ $({\rm mod}\ 4)$. By Lemma~\ref{skla}, we obtain an associated hooked Skolem sequence of order $n$, and by Lemma~\ref{bound}, we obtain a near gracefully labelled Type $(a)$ Dutch windmill with $n$ blocks for $n \geq 22.$ Now, we can see from the above steps and Theorem~\ref{necessarycacti} that a Type $(a)$ Dutch windmill of order $n \geq 22$ with three pendant triangles is graceful or near graceful.
%Following these steps we can show

Skolem and hooked Skolem sequences with three pivots of order $8 \leq n \leq 21$ are given in Table~\ref{ttypea}. By using the same construction as $hS_{7}$ with the given pivots, then we have the labellings for Type $(a)$ Dutch windmill with $n$ blocks for $8 \leq n \leq 21$. In this table and all subsequent tables, we will represent $10$ by $A$, $11$ by $B$, $12$ by $C$, and so on.

Therefore, for $n\geq 6,$ a Type $(a)$ Dutch windmill with $n$ vanes is graceful or near graceful, as required.
\begin{table}[H]
\begin{center}
\scalebox{0.85}{
\begin{tabular}{|r|l|l|}
\hline $n$ & Skolem or hooked Skolem sequence & Pivots \\
\hline
\hline $8$ & $4857411568723263$ & $1,2$ and $4$ \\
\hline $9$ & $759242574869311368$ & $2,4$ and $5$ \\
\hline $10$ & $A853113598A7426249706$ & $5,8$ and $A$ \\
\hline $11$ & $B68527265A8B7941134A309$ & $2,7$ and $B$ \\
\hline $12$ & $A8531135C8A6B9742624C79B$ & $3,6$ and $A$ \\
\hline $13$ & $B97D4262479B6CA8D5311358AC$ & $3,7$ and $B$ \\
\hline $14$ & $CA75311357EACDB864292468EBD09$ & $4,6$ and $C$ \\
\hline $15$ & $FDB9753EC3579BDF6A84CE64118A202$ & $4,9$ and $F$ \\
\hline $16$ & $9FDBG864292468BDFECAG75311357ACE$ & $3,4$ and $F$ \\
\hline $17$ & $FDB9H64282469BDF8GECAH75311357ACEG$ & $4,5$ and $9$ \\
\hline $18$ & $GEC9753113579ICEGHFDA8642B2468AIDFH0B$ & $4,6$ and $9$ \\
\hline $19$ & $IGEC9753113579JCEGIHFDA8642B2468AJDFH0B$ & $7,9$ and $C$ \\
\hline $20$ & $JHFDB9753CK3579BDFHJICEGA86411K468A2E2IG$ & $1,6$ and $9$ \\
\hline $21$ & $KIGEC9753113579LCEGIKAJHFDB8642A2468LBDFHJ$ & $5,A$ and $K$ \\
\hline
\end{tabular}}
\caption{Skolem and hooked Skolem sequences with three pivots for Type $(a)$.}
\label{ttypea}
\end{center}
\end{table}
\vspace{-2cm}
\end{pf}
%%%%%%%%%%%%%%%%%%%%%%%%%%%%%% {{  Type 2 }} %%%%%%%%%%%%%%%%%%%%%%%%%%%%%%%%%%%%%%%%%%%%%%%%%%%%%%%%%%%%%%%%%%%%%%%%%%
\subsection{Type $(b)$}
 For $n \leq 4$, there are not enough triangles to form a Type $(b)$ Dutch windmill.
\begin{lemma}\label{typeb}
Any Type $(b)$ Dutch windmill with at least $5$ blocks is graceful or near graceful.
\end{lemma}
\begin{pf}
For Type $(b)$, consider $S_{5}=\left(4,1,1,5,4,2,3,2,5,3\right)$. This sequence has three pivots, which are $1,2$, and $4$. This gives us a graceful labelling of a Type $(b)$ Dutch windmill with $5$ blocks. The rest of the proof is analogous to the steps taken to prove Lemma~\ref{typea}, by using Lemma~\ref{bound} with $S_{5}$ to obtain a (near) graceful labelling of a Type $(b)$ Dutch windmill with $n$ blocks for $n \geq 16$. Table~\ref{ttypeb} provides Skolem and hooked Skolem sequences with three pivots of order $6 \leq n \leq 15,$ each of which gives a (near) graceful labelling of a Type $(b)$ Dutch windmill with $n$ blocks.
\begin{table}[H]
\begin{center}
\scalebox{0.85}{
\begin{tabular}{|r|l|l|}
\hline $n$ & Skolem or hooked Skolem sequence & Pivots \\
\hline
\hline $6$ & $5611453643202$ & $1,3$ and $5$ \\
\hline $7$ & $746354376511202$ & $1,4$ and $7$ \\
\hline $8$ & $3723258476541186$ & $1,3$ and $5$ \\
\hline $9$ & $759242574869311368$ & $1,4$ and $5$ \\
\hline $10$ & $2529115784A694738630A$ & $2,4$ and $5$ \\
\hline $11$ & $B68527265A8B7941134A309$ & $5,7$ and $B$ \\
\hline $12$ & $A8531135C8A6B9742624C79B$ & $5,6$ and $A$ \\
\hline $13$ & $B97D4262479B6CA8D5311358AC$ & $1,6$ and $7$ \\
\hline $14$ & $CA75311357EACDB864292468EBD09$ & $3,4$ and $C$ \\
\hline $15$ & $FDB9753EC3579BDF6A84CE64118A202$ & $6,9$ and $F$ \\
\hline
\end{tabular}}
\caption{Skolem and hooked Skolem sequences with three pivots for Type $(b)$.}
\label{ttypeb}
\end{center}
\end{table}
\vspace{-2cm}
\end{pf}
%\newline  When considering $n \leq 4$, there are not enough triangles. Now we have shown the associated Skolem and hooked Skolem sequences for all $n$.
%The Skolem sequence $\{(7,8),(3,5),(1,4),(2,6)\}$ has double pivots, $2$ and $3$. The same argument as the independent pendant section applies and we get that Dutch windmills of order at least $13$ with double pendant triangles. Table~\ref{doubleparticular} has Skolem and hooked Skolem sequences with double pivots of order $5$ through $12$. For $n=2$, there are not enough triangles, for $n=3$, the case is the same as a Dutch windmill of order $3$ with no pendant triangles, and for $n=4$, the case is the same as a Dutch windmill of order $3$ with one pendant triangle.

%- give small cases
%- give the base skolem seq
% by theorem number gives us $$a \leq n \leq b$$
%- the remaining cases from table number
%%%%%%%%%%%%%%%%%%%%%%%%%%%%%% {{  Type 3 }} %%%%%%%%%%%%%%%%%%%%%%%%%%%%%%%%%%%%%%%%%%%%%%%%%%%%%%%%%%%%%%%%%%%%%%%%%%
\subsection{Type $(c)$}
 For $n \leq 4$, there are not enough triangles to form a Type $(c)$ Dutch windmills.
\begin{lemma}\label{typec}
Any Type $(c)$ Dutch windmill with at least $5$ blocks is graceful or near graceful.\end{lemma}
\begin{pf}For Type $(c)$, consider $S_{5}=\left(3,5,2,3,2,4,5,1,1,4\right)$. This sequence has three pivots, which are $1,2$, and $3$. This gives us a graceful labelling of a Type $(c)$ Dutch windmill with $5$ blocks. Following the method of Lemma~\ref{typea}, we use $S_{5}$ to obtain a (near) graceful labelling of a Type $(c)$ Dutch windmill with $n$ blocks for $n \geq 16$. Table~\ref{ttypec} provides Skolem and hooked Skolem sequences with three pivots of order $6 \leq n \leq 15,$ each of which gives a (near) graceful labelling of a Type $(c)$ Dutch windmill with $n$ blocks.
%The rest of the proof is analogous to the steps taken to prove Lemma~\ref{typea}, (again by using Lemma~\ref{bound}), we obtain a (near) graceful labelling of a Type $(c)$ dutch windmill with $n$ blocks for $n \geq 16$.

%For Type $(c)$, the smallest Skolem sequence we found was\newline $S_{5}=\left(3,5,2,3,2,4,5,1,1,4\right).$ This sequence has three pivots, which are $1,2$, and $3$. The rest of the proof is analogous to the steps taken to prove Lemma~\ref{typea}, and we get that a gracefully labelled dutch windmill of Type $(c)$ contains $16$ blocks or more. \newline Table~\ref{ttypec} provides Skolem and hooked Skolem sequences with three pivots of order $6 \leq n \leq 15$. When $n \leq 4$, there are not enough triangles.
\begin{table}[H]
\begin{center}
\scalebox{0.85}{
\begin{tabular}{|r|l|l|}
\hline $n$ & Skolem or hooked Skolem sequence & Pivots \\
\hline
\hline $6$ & $6451146523203$ & $1,5$ and $6$ \\
\hline $7$ & $746354376511202$ & $4,6$ and $7$ \\
\hline $8$ & $3723258476541186$ & $1,4$ and $5$ \\
\hline $9$ & $372329687115649854$ & $2,3$ and $7$ \\
\hline $10$ & $36232A768119574A85409$ & $1,2$ and $3$ \\
\hline $11$ & $35232B549A841167B98A607$ & $1,2$ and $3$ \\
\hline $12$ & $3A232C78119AB7685C49654B$ & $1,2$ and $3$ \\
\hline $13$ & $CA8531135D8AC6B9742624D79B$ & $2,4$ and $6$ \\
\hline $14$ & $3B932A211DE9B6CA485647D5E8C07$ & $4,6$ and $A$ \\
\hline $15$ & $ECA8642D2468ACEFB953D735119B70F$ & $4,5$ and $E$ \\
\hline
\end{tabular}}
\caption{Skolem and hooked Skolem sequences with three pivots for Type $(c)$.}
\label{ttypec}
\end{center}
\end{table}
\vspace{-2cm}
\end{pf}
%%%%%%%%%%%%%%%%%%%%%%%%%%%%%% {{  Type 4 }} %%%%%%%%%%%%%%%%%%%%%%%%%%%%%%%%%%%%%%%%%%%%%%%%%%%%%%%%%%%%%%%%%%%%%%%%%%
\subsection{Type $(d)$}
For $n \leq 4$, there are not enough triangles to form a Type $(d)$ Dutch windmill. For $n=5,$ a graceful labelling of the $5$-vane Dutch windmill appears in \cite{rosa}. When $n=6$, see the near graceful labelling in Figure~\ref{taypedn6}.
\begin{figure}[H]
\begin{center}
\scalebox{0.65}{
\begin{tikzpicture}
    \node[shape=circle,draw=black] (19) at (0,0) { $19$};
    \node[shape=circle,draw=black] (0) at (2,0) { $0$};
    \node[shape=circle,draw=black] (17) at (4,0) { $17$};
    \node[shape=circle,draw=black] (1) at (6,0) { $1$};
    \node[shape=circle,draw=black] (16) at (8,0) { $16$} ;
	   \node[shape=circle,draw=black] (6) at (10,0) { $6$};

   % \node[shape=circle,draw=black] (24) at (16,0) {24} ;
		%666666666666666666666666666666666666666666666666666666666666
		    \node[shape=circle,draw=black] (5) at (1,2) { $5$};
		    \node[shape=circle,draw=black] (11) at (3,2) { $11$};
		    \node[shape=circle,draw=black] (13) at (5,2) { $13$};
        \node[shape=circle,draw=black] (3) at (7,2) { $3$};
		    \node[shape=circle,draw=black] (9) at (9,2) { $9$};

		%4444444444444444444444444444444444444444444444444444444444444

    \path [-] (19) edge node[left,above] { $19$} (0);
    \path [-](19) edge node[left,above] { $14\ $} (5);
		\path [-] (5) edge node[left,above] { $\ 5$} (0);

    \path [-] (0) edge node[left,above] { $11\ $} (11);
    \path [-](0) edge node[left,above] { $17$} (17);
		\path [-] (17) edge node[left,above] { $\ 6$} (11);
		
		\path [-] (17) edge node[left,above] { $16$} (1);
    \path [-](17) edge node[left,above] { $4\ $} (13);
		\path [-] (1) edge node[left,above] { $\ 12$} (13);

    \path [-] (1) edge node[left,above] { $2\ $} (3);
    \path [-](1) edge node[left,above] { $15$} (16);
		\path [-] (3) edge node[left,above] { $\ 13$} (16);

	\path [-] (16) edge node[left,above] { $7\ $} (9);
    \path [-](16) edge node[left,above] { $10$} (6);
		\path [-] (6) edge node[left,above] { $\ 3$} (9);
		
		%%%%%%%%%%%%%%%%%%%%%%%%%%%%%%%%%%%%%%%%%%%%%%%%%%%%%%%%%%%%%%%%%%%
 \node[shape=circle,draw=black] (10) at (5,-2) { $10$};
 \node[shape=circle,draw=black] (2) at (7,-2) { $2$};
\path [-] (2) edge node[left,above] { $\ 1$} (1);
    \path [-](1) edge node[left,above] { $9\ $} (10);
		\path [-] (2) edge node[left,above] { $8$} (10);
\end{tikzpicture}
}
\caption{A near gracefully labelled Type $(d)$ Dutch windmill with $6$ blocks.}\label{taypedn6}
\end{center}
\end{figure}
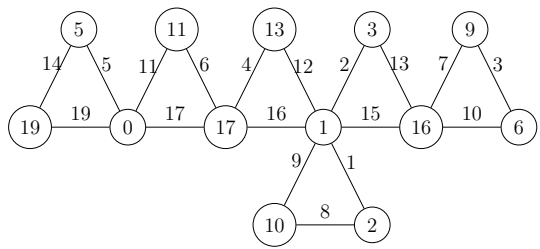
\begin{lemma}\label{typed}
Any Type $(d)$ Dutch windmill with at least $5$ blocks is graceful or near graceful.\end{lemma}
\begin{pf}For Type $(d)$, consider $hS_{7}=\left(3,4,7,3,2,4,2,5,6,7,1,1,5,0,6\right)$. This sequence has three pivots, which are $1,2$, and $3$. This gives us a near graceful labelling of a Type $(d)$ Dutch windmill with $7$ blocks. The rest of the proof is analogous to the steps taken to prove Lemma~\ref{typea}, by using Lemma~\ref{bound} with $hS_{7}$ to obtain a (near) graceful labelling of Type $(d)$ Dutch windmill with $n$ blocks for $n \geq 22$. Table~\ref{ttyped} provides Skolem and hooked Skolem sequences with three pivots of order $8 \leq n \leq 21,$ each of which gives a (near) graceful labelling of a Type $(d)$ Dutch windmill with $n$ blocks.
%For Type $(d)$, the smallest hooked Skolem sequence we found was\newline $hS_{7}=\left(3,4,7,3,2,4,2,5,6,7,1,1,5,0,6\right).$ This sequence has three pivots, which are $1,2$, and $3$. The rest of the proof is analogous to the steps taken to prove Lemma~\ref{typea}, and we get that a gracefully labelled dutch windmill of Type $(d)$ contains $22$ blocks or more. \newline Table~\ref{ttyped} provides Skolem and hooked Skolem sequences with three pivots of order $8 \leq n \leq 21$. When $n \leq 4$, there are not enough triangles, while when $n=5,6$, the graph forms a triangular snake of order $4$, which is graceful [cite ??]. For $n=5$ is gracefully labelled, as can be seen in \cite{rosa}.
%For Type $4$ the smallest hooked Skolem sequence we found is, \newline$hS_{7}=\left(3,4,7,3,2,4,2,5,6,7,1,1,5,0,6\right).$ This sequence has three pivots that are, $1,2$, and $3$. The steps for this part of the proof are analogous to the steps taken the prove Lemma~\ref{type5}.\newline Table~\ref{ttype4} provides Skolem and hooked Skolem sequences with three pivots of order $8 \leq n \leq 21$.
\begin{table}[H]
\begin{center}
\scalebox{0.85}{
\begin{tabular}{|r|l|l|}
\hline $n$ & Skolem or hooked Skolem sequence & Pivots \\
\hline
\hline $8$ & $4857411568723263$ & $1,4$ and $5$ \\
\hline $9$ & $759242574869311368$ & $1,5$ and $7$ \\
\hline $10$ & $A853113598A7426249706$ & $3,8$ and $A$ \\
\hline $11$ & $B68527265A8B7941134A309$ & $2,5$ and $7$ \\
\hline $12$ & $A8531135C8A6B9742624C79B$ & $6,8$ and $A$ \\
\hline $13$ & $B97D4262479B6CA8D5311358AC$ & $4,6$ and $7$ \\
\hline $14$ & $CA75311357EACDB864292468EBD09$ & $1,3$ and $A$ \\
\hline $15$ & $FDB9753EC3579BDF6A84CE64118A202$ & $5,B$ and $F$ \\
\hline $16$ & $9FDBG864292468BDFECAG75311357ACE$ & $4,B$ and $D$ \\
\hline $17$ & $FDB9H64282469BDF8GECAH75311357ACEG$ & $2,6$ and $9$ \\
\hline $18$ & $GEC9753113579ICEGHFDA8642B2468AIDFH0B$ & $3,6$ and $9$ \\
\hline $19$ & $IGEC9753113579JCEGIHFDA8642B2468AJDFH0B$ & $1,5$ and $E$ \\
\hline $20$ & $JHFDB9753CK3579BDFHJICEGA86411K468A2E2IG$ & $6,D$ and $F$ \\
\hline $21$ & $KIGEC9753113579LCEGIKAJHFDB8642A2468LBDFHJ$ & $6,E$ and $G$ \\
\hline
\end{tabular}}
\caption{Skolem and hooked Skolem sequences with three pivots for Type $(d)$.}\label{ttyped}
\end{center}
\end{table}
\vspace{-2cm}
\end{pf}
%\newline When considering $n \leq 4$, there are not enough triangles, while when $n=5$, the graph forms a triangular snake of order $4$, which is graceful [cite ??]. Now we have shown the associated Skolem and hooked Skolem sequences for all $n$.
%%%%%%%%%%%%%%%%%%%%%%%%%%%%%% {{  Type 5 }} %%%%%%%%%%%%%%%%%%%%%%%%%%%%%%%%%%%%%%%%%%%%%%%%%%%%%%%%%%%%%%%%%%%%%%%%%%
\subsection{Type $(e)$}
 For $n \leq 3$, there are not enough triangles to form Type $(e)$ Dutch windmills. For $n=4,$ a graceful labelling of the $4$-vane Dutch windmill appears in \cite{rosa}.
\begin{lemma}\label{typee}
Any Type $(e)$ Dutch windmill with at least $4$ blocks is graceful or near graceful.\end{lemma}

\begin{pf}For Type $(e)$, consider $S_{5}=\left(2,4,2,3,5,4,3,1,1,5\right)$. This sequence has three pivots, which are $1,2$, and $3$. This gives us a graceful labelling of a Type $(e)$ Dutch windmill with $5$ blocks. The rest of the proof is analogous to the steps taken to prove Lemma~\ref{typea}, by using Lemma~\ref{bound} with $S_{5}$ to obtain a (near) graceful labelling of a Type $(e)$ Dutch windmill with $n$ blocks for $n \geq 16$. Table~\ref{ttypee} provides Skolem and hooked Skolem sequences with three pivots of order $6 \leq n \leq 15,$ each of which gives a (near) graceful labelling of a Type $(e)$ Dutch windmill with $n$ blocks.
%For Type $(e)$, the smallest Skolem sequence we found was\newline $S_{5}=\left(2,4,2,3,5,\newline4,3,1,1,5\right).$ This sequence has three pivots, which are $1,2$, and $3$. The rest of the proof is analogous to the steps taken to prove Lemma~\ref{typea}, and we get that a gracefully labelled dutch windmill of Type $(e)$ contains $16$ blocks or more. \newline Table~\ref{ttypee} provides Skolem and hooked Skolem sequences with three pivots of order $6 \leq n \leq 15$. When $n \leq 3$, there are not enough triangles, while when $n=4$, the graph forms a triangular snake of order $4$ with one pendant, which is gracefully labeled, as can be seen in \cite{rosa}.
%For Type $5$ the smallest Skolem sequence we found is, $S_{5}=\left(2,4,2,3,5,\newline4,3,1,1,5\right).$ This sequence has three pivots that are, $1,2$, and $3$. The steps for this part of the proof are analogous to the steps taken the prove Lemma~\ref{type1}.
%\newline Table~\ref{ttype5} provides Skolem and hooked Skolem sequences with three pivots of order $6 \leq n \leq 15$.
\begin{table}[H]
\begin{center}
\scalebox{0.85}{
\begin{tabular}{|r|l|l|}
\hline $n$ & Skolem or hooked Skolem sequence & Pivots \\
\hline
\hline $6$ & $5611453643202$ & $3,4$ and $5$ \\
\hline $7$ & $746354376511202$ & $1,6$ and $7$ \\
\hline $8$ & $3723258476541186$ & $1,2$ and $3$ \\
\hline $9$ & $572825967348364911$ & $2,3$ and $7$ \\
\hline $10$ & $36232A768119574A85409$ & $2,3$ and $6$ \\
\hline $11$ & $B68527265A8B7941134A309$ & $4,6$ and $B$ \\
\hline $12$ & $A8531135C8A6B9742624C79B$ & $4,5$ and $6$ \\
\hline $13$ & $39B32D258A9C5B6784DA647C11$ & $4,8$ and $B$ \\
\hline $14$ & $DB964E1146C9BDA8537E35C8A7202$ & $3,6$ and $B$ \\
\hline $15$ & $CE3693BF262DC97EAB1187F4D5A4805$ & $1,7$ and $B$ \\
\hline
\end{tabular}}
\caption{Skolem and hooked Skolem sequences with three pivots for Type $(e)$.}\label{ttypee}
\end{center}
\end{table}
\vspace{-2cm}
\end{pf}
%%%%%%%%%%%%%%%%%%%%%%%%%%%%%% {{  Type 6 }} %%%%%%%%%%%%%%%%%%%%%%%%%%%%%%%%%%%%%%%%%%%%%%%%%%%%%%%%%%%%%%%%%%%%%%%%%%
\subsection{Type $(f)$}
For $n \leq 3$, there are not enough triangles to form Type $(f)$ Dutch windmills. For $n=4,$ a graceful labelling of the $4$-vane Dutch windmill appears in \cite{rosa}.
\begin{lemma}\label{typef}
Any Type $(f)$ Dutch windmill with at least $4$ blocks is graceful or near graceful.\end{lemma}
\begin{pf}For Type $(f)$, consider $S_{5}=\left(2,4,2,3,5,4,3,1,1,5\right)$. This sequence has three pivots, which are $1,3$, and $4$. This gives us a graceful labelling of a Type $(f)$ Dutch windmill with $5$ blocks. The rest of the proof is analogous to the steps taken to prove Lemma~\ref{typea}, by using Lemma~\ref{bound} with $S_{5}$ to obtain a (near) graceful labelling of a Type $(f)$ Dutch windmill with $n$ blocks for $n \geq 16$. Table~\ref{ttypef} provides Skolem and hooked Skolem sequences with three pivots of order $6 \leq n \leq 15,$ each of which gives a (near) graceful labelling of a Type $(f)$ Dutch windmill with $n$ blocks.
%For Type $(f)$, the smallest Skolem sequence we found was\newline $S_{5}=\left(2,4,2,3,5,\newline4,3,1,1,5\right).$ This sequence has three pivots, which are $1,3$, and $4$. The rest of the proof is analogous to the steps taken to prove Lemma~\ref{typea}, and we get that a gracefully labelled dutch windmill of Type $(f)$ contains $16$ blocks or more. \newline Table~\ref{ttypef} provides Skolem and hooked Skolem sequences with three pivots of order $6 \leq n \leq 15$. When $n \leq 3$, there are not enough triangles, while when $n=4$, the graph forms a triangular snake of order $4$ with two pendant triangles, which is gracefully labeled, as can be seen in \cite{rosa}.
%For Type $6$ the smallest Skolem sequence we found is, $S_{5}=\left(2,4,2,3,5,\newline 4,3,1,1,5\right).$ This sequence has three pivots that are, $1,3$, and $4$. The steps for this part of the proof are analogous to the steps taken the prove Lemma~\ref{type1}.\newline Table~\ref{ttype6} provides Skolem and hooked Skolem sequences with three pivots of order $6 \leq n \leq 15$.
\begin{table}[H]
\begin{center}
\scalebox{0.85}{
\begin{tabular}{|r|l|l|}
\hline $n$ & Skolem or hooked Skolem sequence & Pivots \\
\hline
\hline $6$ & $6451146523203$ & $2,5$ and $6$ \\
\hline $7$ & $746354376511202$ & $5,6$ and $7$ \\
\hline $8$ & $3723258476541186$ & $4,5$ and $7$ \\
\hline $9$ & $746394376825291158$ & $2,6$ and $7$ \\
\hline $10$ & $A869117468A4973523205$ & $2,3$ and $A$ \\
\hline $11$ & $378392A2768B5946A54110B$ & $4,5$ and $9$ \\
\hline $12$ & $2529115B86AC947684B3A73C$ & $4,6$ and $9$ \\
\hline $13$ & $9BD3753AC957B82D2A46C84116$ & $2,7$ and $9$ \\
\hline $14$ & $E7D6C5B47654A8EDCB2928A311309$ & $5,6$ and $7$ \\
\hline $15$ & $3C9382B2E1198CDFAB5647E546AD70F$ & $4,5$ and $B$ \\
\hline
\end{tabular}}
\caption{Skolem and hooked Skolem sequences with three pivots for Type $(f)$.}
\label{ttypef}
\end{center}
\end{table}
\vspace{-2cm}
\end{pf}
%\newline  When considering $n \leq 3$, there are not enough triangles, while when $n=4$, the graph forms a triangular snake of order $4$ with two pendant triangles, which is graceful [cite ??].  Now we have shown the associated Skolem and hooked Skolem sequences for all $n$.
%%%%%%%%%%%%%%%%%%%%%%%%%%%%%% {{  Type 7 }} %%%%%%%%%%%%%%%%%%%%%%%%%%%%%%%%%%%%%%%%%%%%%%%%%%%%%%%%%%%%%%%%%%%%%%%%%%
\subsection{Type $(g)$}
For $n \leq 3$, there are not enough triangles to form Type $(g)$ Dutch windmills. For $n=4,$ a graceful labelling of the $4$-vane Dutch windmill appears in \cite{rosa}.
\begin{lemma}\label{typeg}
Any Type $(g)$ Dutch windmill with at least $4$ blocks is graceful or near graceful.
\end{lemma}
\begin{pf}For Type $(g)$, consider $S_{5}=\left(3,4,5,3,2,4,2,5,1,1\right)$. This sequence has three pivots, which are $2,3$, and $4$. This gives us a graceful labelling of a Type $(g)$ Dutch windmill with $5$ blocks. The rest of the proof is analogous to the steps taken to prove Lemma~\ref{typea},  by using Lemma~\ref{bound} with $S_{5}$ to obtain a (near) graceful labelling of a Type $(g)$ Dutch windmill with $n$ blocks for $n \geq 16$. Table~\ref{ttypeg} provides Skolem and hooked Skolem sequences with three pivots of order $6 \leq n \leq 15,$ each of which gives a (near) graceful labelling of a Type $(g)$ Dutch windmill with $n$ blocks.
%For Type $(g)$, the smallest Skolem sequence we found was\newline $S_{5}=\left(3,4,5,3,2,4,2,5,1,1\right).$ This sequence has three pivots, which are $2,3$, and $4$. The rest of the proof is analogous to the steps taken to prove Lemma~\ref{typea}, and we get that a gracefully labelled dutch windmill of Type $(g)$ contains $16$ blocks or more. \newline Table~\ref{ttypeg} provides Skolem and hooked Skolem sequences with three pivots of order $6 \leq n \leq 15$. When $n \leq 3$, there are not enough triangles, while when $n=4$, the graph forms a triangular snake of order $4$, which is gracefully labeled, as can be seen in \cite{rosa}.
%For Type $7,$ the smallest Skolem sequence we found is, $S_{5}=\left(3,4,5,3,2,4,2,5,1,1\right).$ This sequence has three pivots that are, $2,3$, and $4$. The steps for this part of the proof are analogous to the steps taken the prove Lemma~\ref{type1}.\newline Table~\ref{ttype7} provides Skolem and hooked Skolem sequences with three pivots of order $6 \leq n \leq 15$.
\begin{table}[H]
\begin{center}
\scalebox{0.85}{
\begin{tabular}{|r|l|l|}
\hline $n$ & Skolem or hooked Skolem sequence & Pivots \\
\hline
\hline $6$ & $6451146523203$ & $1,2$ and $4$ \\
\hline $7$ & $746354376511202$ & $1,4$ and $5$ \\
\hline $8$ & $4857411568723263$ & $1,2$ and $5$ \\
\hline $9$ & $759242574869311368$ & $1,4$ and $7$ \\
\hline $10$ & $2529115784A694738630A$ & $1,4$ and $5$ \\
\hline $11$ & $B68527265A8B7941134A309$ & $2,5$ and $B$ \\
\hline $12$ & $A8531135C8A6B9742624C79B$ & $5,6$ and $8$ \\
\hline $13$ & $B97D4262479B6CA8D5311358AC$ & $1,4$ and $6$ \\
\hline $14$ & $CA75311357EACDB864292468EBD09$ & $1,3$ and $C$ \\
\hline $15$ & $FDB9753EC3579BDF6A84CE64118A202$ & $5,6$ and $F$ \\
\hline
\end{tabular}}
\caption{Skolem and hooked Skolem sequences with three pivots for Type $(g)$.}
\label{ttypeg}
\end{center}
\end{table}
\vspace{-2cm}
\end{pf}
%%%%%%%%%%%%%%%%%%%%%%%%%%%%%% {{  Type 8 }} %%%%%%%%%%%%%%%%%%%%%%%%%%%%%%%%%%%%%%%%%%%%%%%%%%%%%%%%%%%%%%%%%%%%%%%%%%
\subsection{Type $(h)$}
 For $n \leq 3$, there are not enough triangles to form Type $(h)$ Dutch windmills. For $n=4,5$ a graceful labelling of the $4$-vane and $5$-vane Dutch windmill appears in \cite{rosa}.
\begin{lemma}\label{typeh}
Any Type $(h)$ Dutch windmill with at least $4$ blocks is graceful or near graceful.
\end{lemma}
\begin{pf}
For Type $(h)$, consider $hS_{6}=\left(4,5,3,6,4,3,5,1,1,6,2,0,2\right)$. This sequence has three pivots, which are $1,3$, and $4$. This gives us a near graceful labelling of a Type $(h)$ Dutch windmill with $6$ blocks. The rest of the proof is analogous to the steps taken to prove Lemma~\ref{typea}, again by using Lemma~\ref{bound} with $hS_{6}$ to obtain a (near) graceful labelling of a Type $(h)$ Dutch windmill with $n$ blocks for $n \geq 19$. Table~\ref{ttypeh} provides Skolem and hooked Skolem sequences with three pivots of order $7 \leq n \leq 18,$ each of which gives a (near) graceful labelling of a Type $(h)$ Dutch windmill with $n$ blocks.
%For Type $(h)$, the smallest Skolem sequence we found was\newline $hS_{6}=\left(4,5,3,6,4,3,5,1,1,6,2,0,2\right)$. This sequence has three pivots, which are $1,3$, and $4$. The rest of the proof is analogous to the steps taken to prove Lemma~\ref{typea}, and we get that a gracefully labelled dutch windmill of Type $(h)$ contains $19$ blocks or more. \newline Table~\ref{ttypeh} provides Skolem and hooked Skolem sequences with three pivots of order $7 \leq n \leq 18$. When $n \leq 3$, there are not enough triangles, while when $n=4,5$, the graph forms a triangular snake of order $4$ and $5$ with one pendant, which is gracefully labeled, as can be seen in \cite{rosa}.
\begin{table}[H]
\begin{center}
\scalebox{0.75}{
\begin{tabular}{|r|l|l|}
\hline $n$ & Skolem or hooked Skolem sequence & Pivots \\
\hline
%\hline $6$ & $4536435116202$ & $1,3$ and $4$ \\
\hline $7$ & $746354376511202$ & $3,4$ and $5$ \\
\hline $8$ & $1157468543763282$ & $3,4$ and $5$ \\
\hline $9$ & $736931176845924258$ & $1,3$ and $7$ \\
\hline $10$ & $5262854A674981137A309$ & $1,4$ and $6$ \\
\hline $11$ & $635A37659B117A842924B08$ & $5,6$ and $7$ \\
\hline $12$ & $52426549B86CA711983B73AC$ & $4,5$ and $6$ \\
\hline $13$ & $86272C568B75D9A11C34B394AD$ & $1,5$ and $8$ \\
\hline $14$ & $2726AD1176C4EBA458D935C3B8E09$ & $4,6$ and $7$ \\
\hline $15$ & $8C3473D48117ACEFB96D52A2659BE0F$ & $2,5$ and $C$ \\
\hline $16$ & $637A3E6FC7D52A2G58BEC9FD48114B9G$ & $2,5$ and $7$ \\
\hline $17$ & $962D2G5649754EFHD78CAGB1138E3FACHB$ & $4,5$ and $D$ \\
\hline $18$ & $5D9485E4F1198BDG6HIAEC6FB7232A3G7CH0I$ & $8,9$ and $E$ \\
%\hline $19$ & $CA74GB94F73AC3J9BIDHG85FE115682D2J6IH0E$ & $3,7$ and $B$ \\
\hline
\end{tabular}}
\caption{Skolem and hooked Skolem sequences with three pivots for Type $(h)$.}
\label{ttypeh}
\end{center}
\end{table}
\vspace{-2cm}
\end{pf}
%%%%%%%%%%%%%%%%%%%%%%%%%%%%%% {{  Type 9 }} %%%%%%%%%%%%%%%%%%%%%%%%%%%%%%%%%%%%%%%%%%%%%%%%%%%%%%%%%%%%%%%%%%%%%%%%%%
\subsection{Type $(i)$}
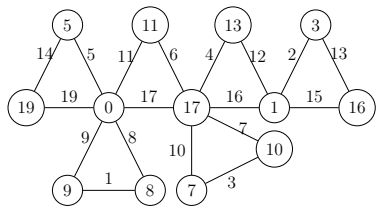
\begin{figure}[H]
\begin{center}
\scalebox{0.55}{
\begin{tikzpicture}
    \node[shape=circle,draw=black] (19) at (0,0) { $19$};
    \node[shape=circle,draw=black] (0) at (2,0) { $0$};
    \node[shape=circle,draw=black] (17) at (4,0) { $17$};
    \node[shape=circle,draw=black] (1) at (6,0) { $1$};
    \node[shape=circle,draw=black] (16) at (8,0) { $16$} ;
	   %\node[shape=circle,draw=black] (6) at (10,0) {\tiny $6$};

   % \node[shape=circle,draw=black] (24) at (16,0) {24} ;
		%666666666666666666666666666666666666666666666666666666666666
		    \node[shape=circle,draw=black] (5) at (1,2) { $5$};
		    \node[shape=circle,draw=black] (11) at (3,2) { $11$};
		    \node[shape=circle,draw=black] (13) at (5,2) { $13$};
        \node[shape=circle,draw=black] (3) at (7,2) { $3$};
		  %  \node[shape=circle,draw=black] (9) at (9,2) {\tiny $9$};

		%4444444444444444444444444444444444444444444444444444444444444

    \path [-] (19) edge node[left,above] { $19$} (0);
    \path [-](19) edge node[left,above] { $14\ $} (5);
		\path [-] (5) edge node[left,above] { $\ 5$} (0);

    \path [-] (0) edge node[left,above] { $11\ $} (11);
    \path [-](0) edge node[left,above] { $17$} (17);
		\path [-] (17) edge node[left,above] { $\ 6$} (11);
		
		\path [-] (17) edge node[left,above] { $16$} (1);
    \path [-](17) edge node[left,above] { $4\ $} (13);
		\path [-] (1) edge node[left,above] { $\ 12$} (13);

    \path [-] (1) edge node[left,above] { $2\ $} (3);
    \path [-](1) edge node[left,above] { $15$} (16);
		\path [-] (3) edge node[left,above] { $\ 13$} (16);

	%\path [-] (16) edge node[left,above] {\tiny $7$} (9);
  %  \path [-](16) edge node[left,above] {\tiny $10$} (6);
	%	\path [-] (6) edge node[left,above] {\tiny $3$} (9);
		
		%%%%%%%%%%%%%%%%%%%%%%%%%%%%%%%%%%%%%%%%%%%%%%%%%%%%%%%%%%%%%%%%%%%
 \node[shape=circle,draw=black] (9) at (1,-2) { $9$};
 \node[shape=circle,draw=black] (8) at (3,-2) { $8$};
\path [-] (8) edge node[left,above] { $1$} (9);
    \path [-](0) edge node[left,above] { $9\ $} (9);
		\path [-] (0) edge node[left,above] { $\ 8$} (8);

 \node[shape=circle,draw=black] (7) at (4,-2) { $7$};
 \node[shape=circle,draw=black] (10) at (6,-1) { $10$};
\path [-] (17) edge node[left] { $\ \ 10$} (7);
    \path [-](7) edge node[left,below] { $3$} (10);
		\path [-] (10) edge node[right] { $7$} (17);

\end{tikzpicture}
}
\caption{A near gracefully labelled Type $(i)$ Dutch windmill with $6$ blocks.}\label{taypein6}
\end{center}
\end{figure}

 For $n \leq 3$, there are not enough triangles to form Type $(i)$ Dutch windmills. For $n=4,5$ a graceful labelling of the $4$-vane and $5$-vane Dutch windmill appears in \cite{rosa}. When $n=6$, see the near graceful labelling in Figure~\ref{taypein6}.
\begin{lemma}\label{typei}
Any Type $(i)$ Dutch windmill with at least $4$ blocks is graceful or near graceful.\end{lemma}
\begin{pf}For Type $(i)$, consider $hS_{7}=\left(7,4,6,3,5,4,3,7,6,5,1,1,2,0,2\right).$ This sequence has three pivots, which are $4,5$, and $6$. This gives us a near graceful labelling of a Type $(i)$ Dutch windmill with $7$ blocks. The rest of the proof is analogous to the steps taken to prove Lemma~\ref{typea}, by using Lemma~\ref{bound} with $hS_{7}$ to obtain a (near) graceful labelling of a Type $(i)$ Dutch windmill with $n$ blocks for $n \geq 22$. Table~\ref{ttypei} provides Skolem and hooked Skolem sequences with three pivots of order $8 \leq n \leq 21,$ each of which gives a (near) graceful labelling of a Type $(i)$ Dutch windmill with $n$ blocks.
\begin{table}[H]
\begin{center}
\scalebox{0.85}{
\begin{tabular}{|r|l|l|}
\hline $n$ & Skolem or hooked Skolem sequence & Pivots \\
\hline
\hline $8$ & $7536835726248114$ & $1,6$ and $7$ \\
\hline $9$ & $975386357946824211$ & $2,6$ and $7$ \\
\hline $10$ & $A853113598A7426249706$ & $2,3$ and $8$ \\
\hline $11$ & $B68527265A8B7941134A309$ & $4,8$ and $B$ \\
\hline $12$ & $A8531135C8A6B9742624C79B$ & $4,6$ and $8$ \\
\hline $13$ & $CA8531135D8AC6B9742624D79B$ & $4,6$ and $A$ \\
\hline $14$ & $DB964E1146C9BDA8537E35C8A7202$ & $3,9$ and $B$ \\
\hline $15$ & $ECA8642D2468ACEFB953D735119B70F$ & $3,A$ and $C$ \\
\hline $16$ & $9FDBG864292468BDFECAG75311357ACE$ & $2,6$ and $9$ \\
\hline $17$ & $FDB9H64282469BDF8GECAH75311357ACEG$ & $3,D$ and $F$ \\
\hline $18$ & $IGECA8642H2468ACEGIFD9753BH357911DF0B$ & $6,7$ and $I$ \\
\hline $19$ & $IGEC9753113579JCEGIHFDA8642B2468AJDFH0B$ & $2,5$ and $E$ \\
\hline $20$ & $JHFDB9753CK3579BDFHJICEGA86411K468A2E2IG$ & $6,F$ and $J$ \\
\hline $21$ & $KIGEC9753113579LCEGIKAJHFDB8642A2468LBDFHJ$ & $6,G$ and $K$ \\
\hline
\end{tabular}}
\caption{Skolem and hooked Skolem sequences with three pivots for Type $(i)$}
\label{ttypei}
\end{center}
\end{table}
\vspace{-1.6cm}
\end{pf}
%%%%%%%%%%%%%%%%%%%%%%%%%%%%%% {{  Type 10 }} %%%%%%%%%%%%%%%%%%%%%%%%%%%%%%%%%%%%%%%%%%%%%%%%%%%%%%%%%%%%%%%%%%%%%%%%%%
\subsection{Type $(j)$ }
 For $n \leq 3$, there are not enough triangles to form Type $(j)$ Dutch windmills. For $n=4$ a graceful labelling of the $4$-vane Dutch windmill appears in \cite{rosa}.
\begin{lemma}\label{typej}
Any Type $(j)$ Dutch windmill with at least $4$ blocks is graceful or near graceful.\end{lemma}
\begin{pf}For Type $(j)$, consider $S_{5}=\left(2,3,2,5,3,4,1,1,5,4\right)$. This sequence has three pivots, which are $1,2$, and $3$. This gives us a graceful labelling of a Type $(j)$ Dutch windmill with $5$ blocks. The rest of the proof is analogous to the steps taken to prove Lemma~\ref{typea},  by using Lemma~\ref{bound} with $S_{5}$ to obtain a (near) graceful labelling of a Type $(j)$ Dutch windmill with $n$ blocks for $n \geq 16$. Table~\ref{ttypej} provides Skolem and hooked Skolem sequences with three pivots of order $6 \leq n \leq 15,$ each of which gives a (near) graceful labelling of a Type $(j)$ Dutch windmill with $n$ blocks.

%For Type $(j)$, the smallest Skolem sequence we found was\newline$S_{5}=\left(2,3,2,5,3,4,1,1,5,4\right).$ This sequence has three pivots, which are $1,2$, and $3$. The rest of the proof is analogous to the steps taken to prove Lemma~\ref{typea}, and we get that a gracefully labelled dutch windmill of Type $(j)$ contains $16$ blocks or more. \newline Table~\ref{ttypej} provides Skolem and hooked Skolem sequences with three pivots of order $6 \leq n \leq 15$. When $n \leq 3$, there are not enough triangles, while when $n=4$, the graph forms a triangular snake of order $4$, which is gracefully labeled, as can be seen in \cite{rosa}.
%For Type $10$ the smallest Skolem sequence we found is, \newline$S_{5}=\left(2,3,2,5,3,4,1,1,5,4\right).$ This sequence has three pivots that are, $1,2$, and $3$. The steps for this part of the proof are analogous to the steps taken the prove Lemma~\ref{type1}.\newline Table~\ref{ttype10} provides Skolem and hooked Skolem sequences with three pivots of order $6 \leq n \leq 15$.
\begin{table}[H]
\begin{center}
\scalebox{0.85}{
\begin{tabular}{|r|l|l|}
\hline $n$ & Skolem or hooked Skolem sequence & Pivots \\
\hline
\hline $6$ & $2326351146504$ & $1,2$ and $3$ \\
\hline $7$ & $232437546115706$ & $2,3$ and $4$ \\
\hline $8$ & $3753811576428246$ & $1,3$ and $5$ \\
\hline $9$ & $759242574869311368$ & $1,2$ and $7$ \\
\hline $10$ & $A853113598A7426249706$ & $1,3$ and $8$ \\
\hline $11$ & $35232B549A841167B98A607$ & $2,4$ and $5$ \\
\hline $12$ & $A8531135C8A6B9742624C79B$ & $3,6$ and $8$ \\
\hline $13$ & $A8D536C358A6B97D42C2479B11$ & $3,4$ and $A$ \\
\hline $14$ & $7A8C53E7358ADB9C6411E469BD202$ & $3,4$ and $A$ \\
\hline $15$ & $DB964F1146E9BD7CA853F735E8AC202$ & $4,5$ and $B$ \\
\hline
\end{tabular}}
\caption{Skolem and hooked Skolem sequences with three pivots for Type $(j)$.}
\label{ttypej}
\end{center}
\end{table}
\vspace{-2cm}
\end{pf}
%%%%%%%%%%%%%%%%%%%%%%%%%%%%%% {{  Type 11 }} %%%%%%%%%%%%%%%%%%%%%%%%%%%%%%%%%%%%%%%%%%%%%%%%%%%%%%%%%%%%%%%%%%%%%%%%%%
\subsection{Type $(k)$}
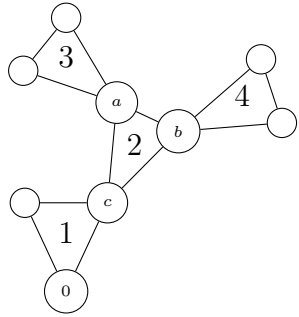
\begin{figure}[H]
\begin{center}
\begin{tikzpicture}[scale=1.3]
\path (0:0cm) node[draw,shape=circle] (v0) {\tiny $0$};
\path (65:1cm) node[draw,shape=circle] (v1) {\tiny $c$};
\path (115:1cm) node[draw,shape=circle] (v2) {$$};
\path (75:2cm) node[draw,shape=circle] (v11) {\tiny $a$};
\path (55:2cm) node[draw,shape=circle] (v12) {\tiny $b$};
\path (101:2.3cm) node[draw,shape=circle] (v3) {$$};
\path (90:2.8cm) node[draw,shape=circle] (v4) {$$};
\path (38:2.8cm) node[draw,shape=circle] (v5) {$$};
\path (50:3.1cm) node[draw,shape=circle] (v6) {$$};

  \node (zz) at (0,0.6) {$1$};
 \node (zz) at (0.7,1.5) {$2$};
 \node (zz) at (0,2.4) {$3$};
 \node (zz) at (1.8,2) {$4$};

%\textsc{1}

\draw
(v0) -- (v1)
(v1) -- (v2)
(v2) -- (v0)
(v1) -- (v11)
(v11) -- (v12)
(v12) -- (v1)
(v11) -- (v3)
(v3) -- (v4)
(v4) -- (v11)
(v12) -- (v5)
(v5) -- (v6)
(v6) -- (v12);
\end{tikzpicture}
\caption{Illustration for Type $(k)$ Dutch windmill labelling.}\label{taypekn}
\end{center}
\end{figure}

Figure~\ref{taypekn} represents a portion of a type $(k)$ Dutch windmill. For this type we can label all the triangles by triples of the form $\left(0,a_{i}+n,b_{i}+n\right)$, where $1 \leq i \leq n,$ using the method of the previous types with the exception of the three pendent triangles $2,3$ and $4$. Let triangle $1$ be labelled by $\left(0,a_{j}+n,b_{j}+n\right)$; we proceed under the assumption that the remaining triangles can be labelled by the pivoting technique. By pivoting we label triangle $2$ by $\left(k,a_{k}+n+k,b_{k}+n+k\right),$ where $b_{j}+n=a_{k}+n+k$ or $a_{j}+n=a_{k}+n+k$. We will consider $b_{j}+n=a_{k}+n+k$ at vertex $c$ in Figure~\ref{taypekn}. Again by pivoting we label triangle $3$ by $\left(l,a_{l}+n+l,b_{l}+n+l\right).$ We will consider $b_{k}+n+k=a_{l}+n+l,$ at vertex $a$, and hence vertex $b$ is labelled $k$. Note that $1 \leq j,k,l \leq n,$ and all are distinct. If we pivot triangle $4$ we label it by the triple $\left(s,a_{s}+n+s,b_{s}+n+s\right).$ However $k \neq s$ and $\min\left(a_{s}+n+s,b_{s}+n+s\right) > n \geq k.$ Since this is impossible we must abandon the pivoting method. Instead we will add some constant $c$ to obtain the triple $\left(c,a_{s}+c+n,b_{s}+c+n\right),$ and we will use triples of the form: $\left(0,i,b_{i}+n\right)$ to avoid any conflicts created between the vertex labels. By Lemma~\ref{lemma2}, adding any constant gives the same differences. (This is a generalization of the idea of pivoting.) The same approach works for hooked Skolem sequences.

For this type, we will introduce the sequences and the corresponding triples to indicate forms of the triples we use.
% In Triangle $4$, we are unable to use the pivot structure, $\left(s,a_{s}+n+s,b_{s}+n+s\right),$ where $1 \leq s \leq n,$ because $s\neq k$ and $n+1 \leq a_{s}+n+s,b_{s}+n+s \leq 3n,$ then instead of using the pivoting structure we will add constant for the triple
%Now, for Triangle $4$ since $1 \leq k \leq n,$ then the pivots structure is replacing the triple $\left(0,a_{i}+n,b_{i}+n\right)$ with $\left(i,a_{i}+i+n,b_{i}+i+n\right),$ but for Type $(k),$ we use $\left(c,a_{i}+c+n,b_{i}+c+n\right),$ for any $c,$ and we could use the both triples forms $\left(0,a_{i}+n,b_{i}+n\right)$ or/and $\left(0,i,b_{i}+n\right).$
%If you see Figure~\ref{taypekn} you can notice that, the vertices $a$ and $c$ can be labelled from the set $\left\{n+1,n+2,\ldots,3n\right\},$ but the vertex $b$ will be labelled from the set $\left\{1,2,\ldots,n\right\},$ this is why we need to use different forms of the triples and the extra shifting technique.
For $n \leq 3$, there are not enough triangles to form Type $(k)$ Dutch windmills. For $n=4$, a graceful labelling of the Dutch windmill with $4$ vanes appears in \cite{rosa}.
\begin{lemma}\label{typek}
Any Type $(k)$ Dutch windmill with at least $5$ blocks is graceful or near graceful.\end{lemma}
\begin{pf}For Type $(k)$, consider $S_{5}=\left(2,3,2,5,3,4,1,1,5,4\right).$ This sequence gives triples of the form $\left(0,a_{i}+n,b_{i}+n\right)$ as follows: $\left(0, 12, 13\right)$, $\left(0, 6, 8\right)$, $\left(0, 7, 10\right)$, $\left(0, 11, 15\right)$, $\left(0, 9, 14\right)$. The above sequence $S_{5}$ has three pivots, but it does not work for this case and so we will use mixed forms of the triples to gracefully label a Type $(k)$ Dutch windmill with $5$ blocks and the triples are: $\left(1,13,14\right),\left(1,7,9\right),\left(3,10,13\right),\left(0,11,15\right),\left(0,5,\right. \\ \left.14\right).$ The first and the third triples are formed by pivoting the elements, the second triple is formed by adding $1$ to each element, the fourth triple is formed by the base block of the form $\left(0,a_{i}+n,b_{i}+n\right)$ and the fifth triple is formed by the base block of the form $\left(0,i,b_{i}+n\right).$ This gives us a graceful labelling of a Type $(k)$ Dutch windmill with $5$ blocks.

We have an order $5$ Skolem sequence $S_{5}$. By Theorem~\ref{langfordexistence}, a hooked Langford sequence of order $m$ with $d=6$ exists for $m \geq 11$. By Lemma~\ref{skla}, we obtain an associated Skolem sequence of order $n,$ and by Lemma~\ref{bound}, we obtain a gracefully labelled Type $(k)$ Dutch windmill with $n$ blocks for $n \geq 16.$ Likewise, By Theorem~\ref{langfordexistence}, a Langford sequence of order $m$ with $d=6$ exists for $m \geq 11$. By Lemma~\ref{skla}, we obtain an associated hooked Skolem sequence of order $n$, and by Lemma~\ref{bound}, we obtain a near graceful labelled Type $(k)$ Dutch windmill with $n$ blocks for $n \geq 16.$ Now, we can see from the above steps and Theorem~\ref{necessarycacti} that a Type $(k)$ Dutch windmill of order $n \geq 16$ with three pendant triangles is graceful or near graceful.
%Following these steps we can show

Skolem and hooked Skolem sequences with corresponding triples of order $6 \leq n \leq 15$ are given in Table~\ref{ttypek}, where the triples in bold are formed in an unusual way (i.e formed by one of the following forms: $\left(i, a_{i}+i+n, b_{i}+i+n\right)$, $(c$, $a_{i}+c+n$, $b_{i}+c+n)$, for some $c,$ or $\left(0,i,b_{i}+n\right)$) and the rest of the triples are formed by the typical $\left(0,a_{i}+n,b_{i}+n\right)$ construction. By using a similar structure as $S_{5}$ with the given triples, then we have the labellings for Type $(k)$ Dutch windmills with $n$ blocks for $6 \leq n \leq 15$.

Therefore, for $n\geq 5,$ a Type $(k)$ Dutch windmill with $n$ vanes is graceful or near graceful, as required.
%The rest of the proof is analogous to the steps taken to prove Lemma~\ref{typea}, (again by using Lemma~\ref{bound}), we obtain a (near) graceful labelling of Type $(k)$ dutch windmill with $n$ blocks for $n \geq 16$. \newline Table~\ref{ttypek} provides Skolem and hooked Skolem sequences with the corresponding triples of order $6 \leq n \leq 15,$ each of which gives a (near) graceful labelling of a Type $(k)$ dutch windmill with $n$ blocks.

%\begin{proof} For Type $(k)$, consider $S_{5}=\left(2,3,2,5,3,4,1,1,5,4\right)$ and its triples as follows:\newline  $\left(1,13,14\right),\left(1,7,9\right),\left(3,10,13\right),\left(0,11,15\right),\left(0,5,14\right).$ These triples label the Dutch windmill of Type $(k)$ of order $5$. \par  The rest of the proof is analogous to the steps taken to prove Lemma~\ref{typea}, and we get that a gracefully labelled dutch windmill of Type $(k)$ contains $16$ blocks or more. Table~\ref{ttypek} provides Skolem and hooked Skolem sequences and the corresponding triples of order $6 \leq n \leq 15$. When $n \leq 3$, there are not enough triangles, while when $n=4$, the graph forms a triangular snake of order $4$ with one pendant, which is which is gracefully labeled, as can be seen in \cite{rosa}.
%\begin{table}[H]

%\begin{center}

%\scalebox{0.75}{
%\begin{tabular}{|r|p{8cm}|p{10cm}|}
%\hline $n$ & Skolem or hooked Skolem sequence & Triples \\
%\hline
\begin{table}[H]

\begin{center}

\scalebox{0.65}{
\begin{tabular}{|r|p{8cm}|p{10cm}|}
\hline $n$ & Skolem or hooked Skolem sequence & Triples \\
\hline

%\hline $5$ & $2325341154$ &  $\left(1,13,14\right),\left(1,7,9\right),\left(3,10,13\right),\left(0,11,15\right),\left(0,5,14\right).$ \\
\hline $6$ & $2326351146504$ & $\left(\textbf{1,14,15}\right),\left(\textbf{1,8,10}\right),\left(\textbf{3,11,14}\right),\left(0,15,19\right),\left(0,12,17\right),\newline \left(\textbf{0,6,16}\right).$  \\

\hline $7$ & $746354376511202$ & $\left(\textbf{0,1,19}\right),\left(\textbf{0,2,22}\right),\left(\textbf{4,15,18}\right),\left(\textbf{7,11,20}\right),\left(0,12,17\right),\newline \left(0,10,16\right), \left(7,15,22\right).$  \\

\hline $8$ & $3753811576428246$ & $\left(\textbf{1,15,16}\right),\left(0,20,22\right),\left(\textbf{2,5,14}\right),\left(0,19,23\right),\left(\textbf{5,16,21}\right), \newline \left(0,18,24\right), \left(0,10,17\right),\left(0,13,21\right).$  \\

\hline $9$ & $736931176845924258$ & $\left(\textbf{7,8,23}\right),\left(\textbf{0,2,25}\right),\left(\textbf{3,14,17}\right),\left(0,20,24\right),\left(0,21,26\right),  \newline \left(0,12,18\right), \left(\textbf{7,17,24}\right),\left(0,19,27\right),\left(0,13,22\right).$  \\

\hline $10$ & $5262854A674981137A309$ & $\left(\textbf{1,25,26}\right),\left(0,12,14\right),\left(0,26,29\right),\left(\textbf{1,5,22}\right),
\left(0,11,16\right), \newline \left(\textbf{6,19,25}\right),\left(0,20,27\right),\left(0,15,23\right),\left(\textbf{0,9,31}\right),\left(0,18,28\right).$  \\

\hline $11$ & $635A37659B117A842924B08$ & $\left(0,22,23\right),\left(0,28,30\right),\left(0,13,16\right),\left(0,27,31\right),
\left(\textbf{7,12,26}\right), \newline \left(\textbf{6,18,24}\right),\left(\textbf{7,24,31}\right),\left(\textbf{0,8,34}\right),\left(0,20,29\right),\left(0,15,25\right), \newline \left(0,21,32\right).$  \\

\hline $12$ & $A8531135C8A6B9742624C79B$ & $\left(0,17,18\right),\left(0,29,31\right),\left(0,16,19\right),\left(0,28,32\right),\left(\textbf{6,11,26}\right), \newline \left(\textbf{6,30,36}\right),\left(0,27,34\right),\left(\textbf{8,22,30}\right),\left(\textbf{0,9,35}\right),\left(0,13,23\right), \newline \left(0,25,36\right),\left(0,21,33\right).$  \\

\hline $13$ &$3113692D286C7B9A5847D54CBA$  & $\left(\textbf{5,6,21}\right),\left(\textbf{2,22,24}\right),\left(0,14,17\right),\left(0,32,36\right),\left(0,30,35\right), \newline \left(\textbf{6,24,30}\right), \left(0,26,33\right),\left(0,23,31\right),\left(0,19,28\right),\left(0,29,39\right), \newline \left(0,27,38\right),\left(0,25,37\right), \left(\textbf{0,13,34}\right).$   \\

\hline $14$ & $B36A349C647BDAE9578C25211D80E$ & $\left(0,38,39\right),\left(\textbf{2,37,39}\right),\left(\textbf{2,5,21}\right),\left(0,20,24\right),\left(0,31,36\right), \newline \left(0,17,23\right),\left(0,25,32\right),\left(0,33,41\right),\left(\textbf{0,9,30}\right),\left(0,18,28\right), \newline \left(\textbf{11,26,37}\right),\left(0,22,34\right),\left(0,27,40\right),\left(0,29,43\right).$  \\

\hline $15$ & $11FB479242DE76B9CF86A35D3E85C0A$ & $\left(\textbf{2,3,19}\right),\left(0,23,25\right),\left(\textbf{3,40,43}\right),\left(0,20,24\right),\left(0,38,43\right), \newline \left(0,29,35\right),\left(0,21,28\right),\left(0,34,42\right),\left(\textbf{9,31,40}\right),\left(0,36,46\right), \newline \left(\textbf{0,11,30}\right),\left(0,32,44\right),\left(0,26,39\right),\left(0,27,41\right),\left(0,18,33\right).$  \\
\hline
\end{tabular}
}
\caption{Skolem and hooked Skolem sequences with the triples for Type $(k)$.}\label{ttypek}
\end{center}
\end{table}
\vspace{-2cm}
\end{pf}
\section{Some Remarks}
%\section{Conclusions and Future Directions}
In this paper, we proved Rosa's conjecture for a new family of triangular cacti: Dutch windmills of any order with three pendant triangles. This result, combined with those of Dyer, et al \cite{dyer}, gives the following theorem.
% We can conclude the results in this paper in the following theorem.
\begin{theorem}
Every Dutch windmill with at most three pendant triangles is graceful or near graceful.
\end{theorem}
\textbf{Question:} Can we gracefully label a Dutch windmill consisting of $m$ triangles and $n$ $4$-cycles, all joined at a common point?

Langford sequences have been classically used to build new Skolem sequences. In this paper, we use this technique to gracefully label Dutch windmills with three pendant triangles. However, in our construction, the Langford sequence implicitly contains a pivot, which when pivoted can gracefully label Dutch windmills with four pendant triangles. Thus we can in many cases label the Dutch windmills with four pendant triangles, for Dutch windmills of a large order.
\begin{theorem} There exists $M$ such that for every Dutch windmill of order $m > M$ with exactly four pendant triangles, where one pendant is attached to a vane containing no other pendants, is graceful or near graceful.\end{theorem}
For example, let $S_{8}=(4,8,5,7,4,1,1,5,6,8,7,2,3,2,6,3)$ be a Skolem sequence and $L_{9}^{17}=(24$, 17, 21, 22, 18, 14, 11, 19, 25, 23, 10, 20, 9, 16, 13, 15, 12, 11, 17, 14, 10, 9, 18, 21, 24, 22, 19, 13, 12, 16, 15, 20, 23, $25)$ be a Langford sequence. Now, if we take the Skolem sequence $S_{8}$ and the Langford sequence $L_{9}^{17},$ then we obtain a new Skolem sequence $S_{25}=(4$, 8, 5, 7, 4, 1, 1, 5, 6, 8, 7, 2, 3, 2, 6, 3, 24, 17, 21, 22, 18, 14, 11, 19, 25, 23, 10, 20, 9, 16, 13, 15, 12, 11, 17, 14, 10, 9, 18, 21, 24, 22, 19, 13, 12, 16, 15, 20, 23, $25)$ of order $25.$ If we take that triples obtained from $S_{25}$ and pivot the triples $1,2,4$ and $11$ we can gracefully label the Dutch windmill of order $25$ with four pendant triangles.
% piviot 1 on 5
% piviot 2 on 3
% piviot 4 on 6
% piviot 11 on 12

%\begin{pf}
%\end{pf}
%\begin{ex}add graph as an example of Theorem~\ref{CBIBD}\end{ex}
%\\

%a) talk here about the relation between skolem and BIBD and add theorem~\ref{CBIBD}\\
%b) Ringel conjecture\\
%c) relation between cyclic decompostion and CBIBD\\
%d) Graceful labelling implies cyclic decomposition\\
%e)CBIBD implies graceful labelling
%Useful results:Perhaps we need to use the following result, since there is a relation between Hamiltonian and graceful. \textbf{ BIBD$(v,k,\lambda)$  $\Longrightarrow$  Hamiltonian, [Hor\'{a}k~and~Rosa, 1988]}
%============================================================================================
%If I have t vans and each van 3 or 4 cycle then I can gracefully labelled.
%Idea combine windmills of C_3 and C_4.
%\textbf{I need to rewrite and the following theorem / related to add langford for any graph we are work on it now}

Furthermore, using the Langford sequences technique to gracefully label triangular cacti, we obtain the following theorem.
\begin{theorem}\label{newtheorem}
Let $G$ be a graph on $m$ edges that can be (near) gracefully labelled. Let $x$ be a vertex of $G$ that obtains the label zero under some (near) graceful labelling.
\begin{enumerate}
\item If $G$ is near gracefully labelled we can obtain a new graceful labelled graph $G^{\ast}$ of size $m+3(l+1)$ by adding an $l+1$ triangular vanes at the vertex $x$.

\item If $G$ is gracefully labelled we can obtain a new graceful labelled graph $G^{\ast}$ of size $m+3l$ edges by adding an $l$ triangular vanes at the vertex $x$.
% by using a Langford sequence with defect $d=m+1.$  (l=length, m= edges)
\item If $G$ is gracefully labelled we can obtain a new near graceful labelled graph $G^{\ast}$ of size $m+3l$ by adding an $l$ triangular vanes at the vertex $x$.
% by using a hooked Langford sequence with defect $d=m+1.$

\item If $G$ is near gracefully labelled we can obtain a new near graceful labelled graph $G^{\ast}$ of size $m+3(l+1)$ by adding an $l+1$ triangular vanes at the vertex $x$.
%By identify that 0 is center of Dutch windmill with any vertex of $G$ can be labelled zero in some graceful labelling.
% by using an \textbf{(we have to use the same name ?)} extended Langford sequence with defect $d=m+1.$ (\textbf{I have to explain the technique} )
\end{enumerate}
 \end{theorem}
 In \cite {linek} Linek and Jiang studied $p$-extended Langford sequences. Here we will present the definition of a $p$-extended Langford sequence because we will use it in the proof of Theorem~\ref{newtheorem}.

A \textit{$p$-extended Langford sequence} of defect $d$ and $m$ differences is a sequence $S=\left(s_{1},s_{2},\ldots,s_{2m+1}\right)$ which satisfies these conditions:
\begin{enumerate}
\item for every $k\in\left\{d,d+1,\ldots,d+m-1\right\},$ there exist exactly two elements $s_{i},s_{j}\in S$ such that $s_{i}=s_{j}=k$;

\item if $s_{i}=s_{j}=k,$ with $i < j,$ then $j-i=k$;

\item $s_{p}=0$, for some $1 \leq p \leq 2m+1$.
\end{enumerate}
%. A $k-$extended Langford sequence of defect $d$ and $m$ differences is an integer sequence is a sequence $s_{1},s_{2},\ldots,s_{2m+1}$ that satisfies conditions (1) and (2) for a Langford sequence, and in addition satisfies (3) sk = 0 .

%A $k-$extended Langford sequence of defect $d$ and length $m$ is a sequence $s_{1},s_{2},\ldots,s_{2m+1}$ in which $s_{k}=\epsilon$, where $\epsilon$ is the null symbol, and each other member of the sequence comes from the set $\left\{d,d+1,\ldots,d+m-1\right\}$. Each $j \in \left\{d,d+1,\ldots,d+m-1\right\}$ occurs exactly twice in the sequence, and the two occurrences are separated by exactly $j-1$ symbols.
Now we will present the proof of Theorem~\ref{newtheorem}.\\
\begin{pf}
We present a proof for the first case. The proof of the other three cases will be similar.

Let $G$ be a near gracefully labelled graph whose vertex labels are a subset of $\{0$, $1$, $2$, $\ldots$, $m+1\}$ and whose edge labels are exactly $\{1$, 2, $\ldots$, $m-1$, $m+1\}$.

Form a $p$-extended Langford sequence ($ELS$) with defect $d=m+2$ and the differences $d,d+1,\ldots,d+l-1$ and the extended position at the place $p=s_{2l+1}-(m-1)$. Following the constructions in \cite{linek1}, form a $p$-extended Langford sequence with $d=m+2$. Form a new sequence $ELS^{\ast}$ with all of the elements from $ELS$ along with the element $m$ at the position $k$ and $s_{2l+2}.$ Then from $ELS^{\ast}$ we can construct triples of the form $D=\{0,a_{i}+d+l-1,b_{i}+d+l-1\}, 1 \leq i \leq d+l-1$. These base blocks give the vertex labels for $l+1$ triangles, where 0 is repeated $l+1$ times, as a common vertex.
We obtain the vertex labels from the set $\left\{0,d+l,d+l+1,\ldots,d+3l-1\right\},$ where $0$ is a common vertex and edge labels exactly are $\left\{d+l,d+l+1,\ldots,d+3l-1\right\} \cup \left\{m\right\}.$

The base blocks formed by the $ELS^{\ast}$ sequence with the original labelling of the graph $G$ give a graceful labelling for a new graph $G^{\ast}$, formed by attaching $l+1$ triangles at $0$ on the graph $G$.

%For 2, we can use the same technique of the proof, with Langford sequence, for 3, use hooked Langford sequence and for 4, hooked $ELS^{\ast}$.
For the second statement, we use the same technique as in the first with $L$ as a Langford sequence with defect $d=m+1$ and order $l$. For the third statement, we use the same technique as in the first with $hL$ as a hooked Langford sequence with defect $d=m+1$ and order $l$. For the fourth statement, we use the same technique as in the first with $W$ as a hooked extended Langford sequence with defect $d=m+2$ and order $l$ and use $W^{\ast}$.
\end{pf}

In fact, by using Theorem~\ref{newtheorem}, we can (near) gracefully label new graphs. For instance, let $G=C_{5}$ be a nearly gracefully labelled graph where the edge labels are exactly $\left\{1,2,3,4,6\right\}$; see Figure~\ref{c}. Let $ELS_{7}^{13}=(11$, 14, 15, 16, 17, 18, 19, 7, 8, 9, 10, 11, 12, 13, 7, 14, 8, 15, 9, 16, 10, 17, 0, 18, 12, 19, $13)$ be an extended Langford sequence and $(ELS_{7}^{13})^{\ast} = (11$, 14, 15, 16, 17, 18, 19, 7, 8, 9, 10, 11, 12, 13, 7, 14, 8, 15, 9, 16, 10, 17, $\textbf{5}$, 18, 12, 19, 13, $\textbf{5})$ be a modified extended Langford sequence.

This sequence, $(ELS_{7}^{13})^{\ast}$, gives triples of the form $\left(0,a_{i}+19, b_{i}+19\right)$ as follows: $\left(0, 42, 47\right)$, $\left(0, 27, 34\right)$, $\left(0, 28, 36\right)$, $\left(0, 29, 38\right)$, $\left(0, 30, 40\right)$, $\left(0, 20, 31\right)$, $\left(0, 32, 44\right)$, $\left(0, 33, 46\right)$, $\left(0, 21, 35\right)$, $\left(0, 22, 37\right)$, $\left(0, 23, 39\right)$, $\left(0, 24, 41\right)$, $\left(0, 25, 43\right)$, $\left(0, 26, 45\right)$.
The base blocks formed by the $(ELS_{7}^{13})^{\ast}$ sequence with the original labelling of the graph $G=C_{5}$ give a graceful labelling for a new graph $G^{\ast}$, formed by attaching $14$ triangles at $0$ on the graph $G$.

In 1989, Moulton \cite{mou} proved Rosa's conjecture for a triangular snake, a type of triangular cactus whose block cutpoint graph is a path. Finally we pose the following question.

\textbf{ Question:} Can we use Skolem type sequences to gracefully label triangular snakes? Can we use them to gracefully label triangular snakes with pendant triangles?
\section*{Acknowledgement}
Author Danny Dyer acknowledges research grant support from NSERC.

%\fi

\end{document}